\documentclass[12pt,reqno,a4paper]{amsart}
\usepackage{amsmath,amssymb,amsfonts}
\usepackage[mathscr]{eucal}

\newtheorem{proposition}{Proposition}
\newtheorem{theorem}{Theorem}

\newtheorem{lemma}{Lemma}

\theoremstyle{definition}
\newtheorem{remark}{Remark}
\newtheorem{example}{Example}[section]

\newcommand{\cal}{\EuScript}

\renewcommand{\leq}{\leqslant}
\renewcommand{\geq}{\geqslant}

\renewcommand{\L}{{\cal L}}

\newcommand{\RR}{\mathbb R}

\newcommand{\p}{\partial}

\def\a{\alpha}

\def\d{\delta}
\def\t{\tilde}
\def\vare{\varepsilon}
\def\s{\sigma}
\def\spr {\Sigma_{\rm princ.}}

\def\G {\Gamma}
\def\hM {\widehat M}
\def\GG {{\cal G}}
\def\hL {{\widehat L}}
\def\LL {{\cal L}}
\def\hLL {{\widehat{\cal L}}}
\def\X  {{\bf X}}
\def\hX  {\widehat{\bf X}}
\def\hQ {{\hat Q}}
\def\Y  {{\bf Y}}
\def\K  {{\bf K}}

\def\F {{\cal F}}
\def\D {{\cal D}}
\def\V {{\cal V}}
\def\S {{\mathbf S}}
\def\hS {{\widehat S}}
\def\K {{\mathbf K}}
\def\m{\medskip}

\def\l {{\lambda}}
\def\hw{{\widehat w}}
\def\hl{{\widehat\lambda}}

\def\hD{{\widehat {\Delta}}}
\def\wD #1 {{\widehat {\Delta_{#1}}}}
\def\hPi {{\widehat {\Pi}}}
\def\hP {{\widehat {P}}}
\def\hL {{\widehat {L}}}

\def\ad {{\rm ad}}

\def\Diff {{\rm Diff}}

\def\proj {{\rm proj\,}}
\def\Diff {{\rm Diff\,}}
\def\diff {{\rm diff\,}}
\def\SDiff {{\rm SDiff\,}}
\def\sdiff {{\rm sdiff\,}}

\newcommand{\bs}{{\boldsymbol{s}}}

\newcommand{\rh}{{\boldsymbol{\rho}}}

\title[Operator pencils on algebra of...]
{Operator pencils on the algebra of densities}

\author{A.~Biggs}
\author{H.~M.~Khudaverdian}

\address{School of Mathematics,  University of Manchester,
 Oxford Road,  Manchester   M13 9PL,  UK}
\email{khudian@manchester.ac.uk\\adam.biggs@student.manchester.ac.uk }



\begin{document}

\maketitle
\begin{abstract} In this paper we continue to study
equivariant pencil liftings and differential operators
 on the algebra of
densities. We emphasize the role that the geometry of the extended
manifold plays where the algebra of densities is a special class
of functions. Firstly we consider basic examples. We give
a projective line 
of $\diff(M)$-equivariant pencil liftings for first order operators,
and the canonical second order self-adjoint lifting. Secondly we study
 pencil liftings equivariant with respect to volume preserving 
transformations. This helps to understand the role of self-adjointness
for the canonical pencils. Then we introduce the Duval-Lecomte-Ovsienko
(DLO)-pencil lifting which is derived from 
the full symbol calculus of projective
quantisation. We use the DLO-pencil lifting to describe all regular 
$\proj$-equivariant pencil liftings. In particular the 
comparison of these pencils
with the canonical pencil for second order operators leads to objects related
to the Schwarzian.   
\end{abstract}

 \section {Introduction}

  We consider pencils of operators on the algebra of densities
   passing through a given operator. 
   
Firstly let us fix notations.
 Let $\F_\l(M)$ be the space of densities of weight $\l$ on a 
manifold $M$
for $\l$ an arbitrary real number.
  If $\bs=s(x)|Dx|^\l\in \F(M)$ then under changing of 
local coordinates $x=x(x')$,
  \begin{equation*}\label{transformationequations1}
                          \bs
           =s(x)|Dx|^\l=
s\left(x\left(x^\prime\right)\right)\left|\det
 \left({\p x\over \p x'}\right)\right|^\l
                    |D{x^\prime}|^\l\,.
        \end{equation*}
The space of functions on $M$ is $\F_0(M)$,
densities of weight $\l=0$.
We can multiply densities of different weights. If $\bs_1=s_1(x)|Dx|^{\l_1}$
and $\bs_2=s_2(x)|Dx|^{\l_2}$ are densities of weights $\l_1,\l_2$
respectively then  
   $
\bs=\bs_1\cdot\bs_2=s_1(x)s_2(x)|Dx|^{\l_1+\l_2}$
is a density of weight $\l_1+\l_2$.
We come to the algebra $\F(M)=\oplus_\l\cal F_\l(M)$ 
of densities of all weights on the manifold $M$.

 We also consider 
differential operators on the whole algebra of densities.
 The operators of order zero are
operators of multiplication by an arbitrary non-zero density; 
$L$ is an operator of order $\leq n+1$ on $\F(M)$ 
if for an arbitrary 
density $\bs$, the commutator  $[L,\bs]=L\circ \bs-\bs\circ L$ 
is an operator of order $\leq n$.  For example 
consider the linear operator $\hw$ on $\F(M)$ 
defined by
      \begin{equation}\label{weightoperator}
\hw (\bs)=\l\bs\,\quad {\rm if}\,\bs\in\F_\l(M)\,,\quad
   \hbox{$\hw$ is the weight operator}\,.
     \end{equation} 
It is easy to see that $\hw$ is a first order operator
on the algebra of densities, since 
$\hw(\bs_1\cdot\bs_2)=$ $\bs_1\hw(\bs_2)+\bs_2\hw (\bs_1)$.
For example the operator
  $\hD=(\hw^2+1){\p^2\over \p x^2}+{\p\over \p x}$
is a fourth order operator on the algebra of densities on  the line $\RR$.
Its restriction to the space  of densities of weight $\l$
is the operator
        $$
  \Delta_{\l}=\hD\big\vert_{\hw=\l}=
  (\l^2+1){d^2\over dx^2}+{d\over dx}\,,
        $$
which is an operator of order $2$ on the space $\F_\l(\RR)$.
   We denote by $\D^{(n)}(M)$ the space of linear differential operators
 of order $\leq n$ on the space $\F_\l(M)$ of densities
of weight $\l$, and we denote by $\D^{(n)}(\hM)$ 
the space of operators of order $\leq n$ on the algebra $\F(M)$
(the latter notation will become clear soon).  

One can consider a pencil of operators (operator pencil)
$\{\Delta_\l\}$, i.e. a family of operators depending on
the  parameter $\l$, such that for each $\l$, $\Delta_\l$ is a 
differential operator on $\F_\l$. 
An  operator $\hD$ on the algebra of densities $\F(M)$ 
defines the pencil of operators:
          \begin{equation*}
\hD\mapsto \,\,\{\Delta_\l\}\colon \,\,
\Delta_\l=\hD\big\vert_{\hw=\l}\,,
      \end{equation*} 
such that dependence on $\l$ is polynomial.
We consider only pencils depending polynomially on $\l$
and we identify operators $\hD$ on the algebra of densities 
with their corresponding pencils.
(For details see \cite{KhBiggs1} and \cite{KhVor5}.)

It is useful to
consider the fibre bundle $\hM\to M$ which is the frame bundle of 
the determinant bundle over the manifold $M$.
One can identify the algebra of densities $\F(M)$ with  
the algebra of quasipolynomials  on 
the extended manifold $\hM$ (see \cite{KhVor2}). 
 In local coordinates use a formal variable
$t$ instead $|Dx|$. An arbitrary density 
$\bs\in \F(M)$ can be identified 
with function
$\sum s_r(x)t^{\l_r}$ on the extended manifold $\hM$,
which is quasipolynomial on fibre variable $t$.
  We use the local coordinates $(x^i, t)$ 
on $\hM$.
 Local coordinates in the bundle change in the following way:
     \begin{equation}\label{changingoflocalcoordinatesinextendedmanifold}
(x^{i'},t')\colon\quad
   x^{i'}=x^{i'}(x^i),\quad t'=t'(x^{i},t)=
\det\left({\p x^{i'}\over \p x^{i}}\right)t\,.        
\end{equation}
Differential operators on the algebra of densities can be
identified with differential operators
on quasipolynomials. In particular the weight operator,
 \eqref{weightoperator}, has the form
$\hw=t{\p\over \p t}$. 

 Differential operators on algebra of densities
have natural grading
by their weights. We say that an operator $\hD$ has 
weight $\d$ if it maps the spaces $\F_\l(M)$ to the spaces $\F_{\l+\d}(M)$.
(An operator $\hD$ has weight $\d$ iff 
$[\hw, \hD ]=\hw\circ\hD-\hD\circ\hw=\d\hD$.)

It is also important to distinguish  {\it vertical operators}. 
An operator $\hD$ is a vertical operator if it 
commutes with multiplication  
by an arbitrary function $f$: $\hD(f\bs)=f\hD(\bs)$.  
A differential operator $\hD$ of weight $\d$ 
has the local appearance $\hD=t^\d\hD^\prime\left(x^k,
      \p_i,\hw\right)$.
We will focus on operators of weight $\d=0$.
 Vertical operator of weight $\d=0$ have the global 
appearance $\hD=\sum\hw^ic_i(x)$,
where the $c_i(x)$ are scalar functions. 

    The algebra $\F(M)$ of densities is endowed with a 
canonical scalar product.

If  $\bs_1=s_1(x)|D{x}|^{\l_1}$ and $\bs_2=s_2(x)|D{x}|^{\l_2}$ 
are two densities with compact support then 
                \begin{equation}\label{canonicalscalarproduct}
                   \langle \bs_1,\bs_2\rangle
                          =
                      \begin{cases}
                        & \int_M s_1(x)s_2(x)|D{x}|\,,
                 \quad {\rm if}\quad\l_1+\l_2=1\,,\cr
                                                               \cr
                        &  0\qquad 
          {\rm if}\qquad \l_1+\l_2\not=1\,.\cr
                      \end{cases}
                 \end{equation}
This construction  turns out to be a very 
important tool when studying the geometry 
of differential
operators on $M$ (see \cite{KhVor2}, \cite{KhVor4}.)   
In particular  the scalar
product defines the adjoint of 
linear operators on the algebra of densities. 
 A linear operator $\hD$ acting on densities has an adjoint 
$\hD^*\colon$ $\langle\hD\bs_1,\bs_2\rangle=
\langle\bs_1,\hD^*\bs_2\rangle$. It is evident
that $(x^k)^*=x^k$ and $\p_i^*=-\p_i$.
  One can see that for the weight operator $\hw$, $\hw^*=1-\hw$.
Indeed let $\bs_1$ and $\bs_2$ be
two arbitrary densities of weights $\l_1$
and $\l_2$ respectively. Then
$\langle\hw\bs_1,\bs_2\rangle=\l_1\langle\bs_1,\bs_2\rangle$,
and $\langle\bs_1,\hw^*\bs_2\rangle=$
$\langle\bs_1,(1-\hw)\bs_2\rangle=(1-\l_2)\langle\bs_1,\bs_2\rangle$,
and $\l_1\langle\bs_1,\bs_2\rangle=(1-\l_2)\langle\bs_1,\bs_2\rangle$.
(If $\l_1+\l_2\not=1$ then 
$\langle\bs_1,\bs_2\rangle=0$, if $\l_1+\l_2=1$ then $\l_1=1-\l_2$.)	
\footnote{Note that bilinear form (3) is not positive definite. E.g. for
density $\bs=f(x){|Dx|^\l-|Dx|^{1-\l}\over \sqrt 2}$, 
($\l\not={1\over 2}$) we have that 
$\langle\bs,\bs\rangle=-\int_M f^2(x)Dx$.}

  Our task is to construct a map $\hPi$ 
   which maps   
each operator on densities of a given weight $\l_0$ 
to a pencil of operators which
passes through this operator:
      \begin{equation}\label{pencilpassingthroughoperator}
\hPi\colon \quad  
 \D_{\l_0}^{(n)}\ni \Delta\mapsto \hD=
  \hPi(\Delta)=\{\Delta_\l\}\,,
\quad\hbox {such that}\quad
 \hPi(\Delta)\big\vert_{\hw=\l_0}=\Delta\,.
      \end{equation}
We call such a map $\hPi$ a pencil lifting map.

 We say that a pencil lifting map $\hPi$ is {\it regular} on 
the space of operators of order $\leq n$ if
this map  takes values in operators of order $\leq n$:
    \begin{equation*}
   \hPi\colon\quad  \D^{(n)}_{\l_0} (M)\to\D^{(n)} (\hM)\,.
     \end{equation*} 
  To make the problem well-defined we put restrictions on
the pencil lifting map.  We consider {\it regular} pencil lifting maps
which are {\it equivariant with respect to some Lie subalgebra 
  of vector fields} \cite{KhBiggs1}.  
  Let us briefly recall what this means.

  A pencil lifting map is 
 equivariant with respect to a subgroup $G$ of the group 
$\Diff(M)$ of diffeomorphisms of $M$, i.e.  it is a $G$-lifting,  
if for every diffeomorphism $\varphi\in G$, $\hPi(\Delta^\varphi)=
\left(\hPi(\Delta)\right)^\varphi$. If $\varphi$ is a diffeomorphism which
is infinitesimally closed to the identity, $\varphi=1+\vare\K$,
where $\K$ is a vector field on $M$, then invariance with 
respect to diffeomorphism
$\varphi$ implies that
      \begin{equation}\label{vectorfieldinvariant}
    \ad_\K\left(\hPi(\Delta)\right)=\hPi\left(\ad_\K(\Delta)\right)\,,
     \end{equation}
where operation $\ad_\K$ is defined by condition that
for every density $\bs$, 
 $\,\ad_\K\left(\Delta\right)\bs=$
   $\L_\K(\Delta\bs)-\Delta\left(\L_\K\bs\right)$.
($\L_\K$ is the Lie derivative with respect to vector field $\K$.)
  We say that a pencil lifting map $\hPi$ is 
equivariant with respect to  
Lie algebra $\GG\subseteq \diff(M)$, a $\GG-lifting$,
if the condition \eqref{vectorfieldinvariant} holds
for every vector field $\K\in\GG$.

\begin{remark} 
 We consider regular pencil lifting maps 
which are equivariant with respect to Lie algebras. 
 Often the results will be valid for corresponding groups.
\end{remark}

Here we consider pencil lifting maps 
which are equivariant with respect
to 

\begin{itemize}

\begin{item}
  the algebra $\diff(M)$ of all vector fields on $M$,

\end{item}
\begin{item}
  the algebra $\sdiff(M)=\sdiff_\rh(M)$ of all divergenceless vector fields
   on $M$ with respect to some volume form $\rh$,

\end{item}
\begin{item}
  the algebra $\proj(M)$ of projective vector 
fields on a projective manifold $M$. 

\end{item}

\end{itemize}
 Recall that
  a vector field $\X$ is divergenceless with respect to
the volume form $\rh$, ($\rh\in \F_1(M)$) if the volume form $\rh$ 
is invariant with respect to the vector field $\X$, $\LL_\X\rh=0$.
The divergence of a vector field with respect to 
$\rh$ 
 is defined by the formula
      \begin{equation*}
   {\rm div\,}_\rh \X={\LL_\X\rh\over \rh}=\p_iX^i(x)+
    X^i(x)\p_i\log\rho(x)\,.
     \end{equation*}
The Lie algebra $\sdiff_\rh(M)$ is a Lie subalgebra 
of all vector fields 
of $M$ such that they preserve the volume form $\rh$ on $M$.

A projective  
manifold of dimension $d$ which is locally a projective  
space $\RR P^d$. 
In local projective coordinates $x^i=(x^1,\dots,x^d)$
the Lie algebra $\proj(M)$ is generated by the vector fields 
     \begin{equation*}
\p_i, \qquad x^i\p_k, \qquad     x^ix^k\p_k\,,\qquad
(i,k=1,\dots,d)\,.
      \end{equation*}
(Two different local projective coordinates 
in a vicinity of a point are related by equations
 $x^{i'}={a^i_jx^j+b^i\over c_kx^k+d}$, $i,i'=1,\dots.d$.)

   In what follows we consider 
regular pencil liftings equivariant with respect
to these algebras. 
We also focus our attention on the self-adjointness of
pencil lifting maps.

It should be noted that regular pencil lifting maps
can sometimes be described by a full symbol map.
Let us explain what this means in more detail.

  Let $\GG$ be an arbitrary 
Lie subalgebra of vector fields on $M$. 
  Consider (if it exists) a $\GG$-equivariant {\it full symbol map} 
 $\s^{(n)}_\l$ from the space $\D^{(n)}_\l(M)$
of differential operators (of weight $\d=0$)  on densities of weight $\l$
to the symbol space of contravariant symmetric
  tensor fields of rank $\leq n$ on $M$. 
     The symbol space can be identified
with the space of functions on $T^*M$ which are 
polynomials of order $\leq n$ in the fibre  variables (momenta).
The full symbol map $\s_\l^{(n)}$ is a $\GG$-equivariant map
that prolongs the principal symbol map:
     \begin{equation*}\label{algebra4invariance}        
\hbox {for arbitrary}\,\,\K\in\GG,\quad
  \s^{(n)} (\ad_\K\Delta)=\L_\K\s^{(n)}_\l(\Delta)\,,      
\end{equation*}
and
          \begin{equation}\label{fullsymbolmap}
\hbox {for arbitrary}\,\,
 \Delta\in \D^{(k)}_\l(M)\subseteq \D_\l^{(n)}(M) (k\leq n)\,,\,
 \s^{(n)}_\l(\Delta)=\spr^{(k)}(\Delta)+\dots\,,
          \end{equation}
where the ellipsis means tensors of rank $< k$, and $\spr^{(k)}$
is the principal symbol: for arbitrary
 $\Delta=S^{i_1\dots i_k}\p_{i_1}\dots\p_{i_k}+\dots $,
               \begin{equation*}\label{principalsymbol}
   \spr^{(k)}(\Delta)=S^{i_1\dots i_k}\xi_{i_1}\dots\xi_{i_k}\,.
          \end{equation*}
  The inverse of such a full symbol map defines
a $\GG$-equivariant 
quantisation map  $Q^{(n)}_\l=\left(\s^{(n)}_\l\right)^{-1}$ 
on the space of symbols. (See the analysis in 
section \ref{projective} for details). 

  It is well-known that there is no $\diff(M)$-equivariant full symbol
map for $n\geq 2$. On the other hand there exists a full symbol map
$\sigma^{(n)}_\l$ for arbitrary $n$ and $\l$ in the case if 
$\GG$ is the algebra of projective vector fields, or 
 the algebra of conformal vector 
 fields on $\RR^d$, and this full symbol map is uniquely 
defined by its equivariance with respect to this algebra.
 These are famous results on equivariant quantisation
 due to Duval, Lecomte and Ovsienko 
(see \cite{DuvLecomteOvs} and \cite{LecomteOvs1} 
and also the book \cite{OvsTab}).

The existence of $\GG$-equivariant full symbol map $\sigma^{(n)}_\l$
for all $\l$, determines 
 the corresponding 
pencil lifting map on the spaces $\D^{(n)}_\l(M)$ of operators on densities
of weight $\l$:
       \begin{equation}\label{DLOpencillifting1}
 \hPi_\l\colon\,\, \forall\Delta\in \D^{(n)}_\l(M)\,,\quad
          \left(\hPi_\l(\Delta)\right)\big\vert_{\hw=\mu}= 
      Q^{(n)}_\mu\circ \s^{(n)}_{\l}(\Delta)\,,\quad
        Q^{(n)}_\mu=(\s^{(n)}_{\mu})^{-1}\,.
       \end{equation}
This pencil lifting can be defined for arbitrary $n$. It is evidently
equivariant with respect to the algebra $\GG$. 
One can see that
this is regular pencil lifting. We call this pencil lifting
the Duval-Lecomte-Ovsienko (DLO)-pencil lifting (see in more details
 in the section \ref{projective}).

The DLO-pencil lifting
can be factored through a symbol map.
Of course not every pencil lifting can be factored through 
symbol  map. 
  Recall that regular 
pencil $\hPi$ lifting
defined on $\D^{(n)}_\l(M)$ is {\it strictly regular}
(see \cite{KhBiggs1})
 if its restriction on any subspace 
$\D^{(k)}_\l(M)\subseteq \D^{(n)}_\l(M)$ is also regular:
         \begin{equation}\label{strictlyregular}
\hPi\colon\quad \forall k\leq n\,,
        \D^{(k)}_\l(M)\to \D^{(k)}_\l(\hM)
         \end{equation}
  One can see that DLO-pencil pencil lifting is not only regular but
also strictly regular pencil lifting on all the spaces $\D^{(n)}_\l(M)$.

 One can see that a regular pencil lifting which can be factored through
a symbol map is strictly regular. Being strictly regular 
 is a necessary but not a sufficient condition for a pencil lifting 
to be factored. 
However this simple criterion allows us to 
observe the fact that a pencil lifting cannot be factored through the
 full symbol map. 
\smallskip

   The paper is set out in the following way:
  
In the next two sections   
we consider regular $\diff(M)$-pencil liftings of operators of order $n=1,2$.
Firstly  we consider regular
pencil liftings for first order operators and come to a projective
line of $\diff(M)$-equivariant liftings.
 Next we study the canonical pencil for second order 
operators. Here we construct this canonical pencil
as a Laplace-Beltrami like-operator on extended manifold
$\hM$ and study its geometry
using a Kaluza-Klein like formalism.  

 Higher order regular $\diff(M)$-pencil liftings ($n\geq 3$) 
are impossible  in general 
(see \cite{LecomteMath1}, \cite{Math2}).       
In the fourth section we describe the family of all
regular pencil liftings equivariant with 
respect to volume preserving transformations.
Unlike $\diff(M)$-pencil liftings these pencil liftings
are defined for operators of all orders. 
    Within this family is a
distinguished lifting which can be used to study the role of 
self-adjointness for canonical pencils.

  In the fifth section using  the projective full symbol map
we study the Duval-Lecomte-Ovsienko (DLO)-pencil lifting, 
 and on the base
 of it we construct all regular self-adjoint and anti-self-adjoint
$\proj$-equivariant  pencil liftings.  
In particular we analyze a family of self-adjoint pencil
liftings for second order operators, and comparing it with the canonical
self-adjoint pencil lifting come to geometrical objects whose
transformation law is governed by a multidimensional Schwarzian.

   Throughout this paper by default  we consider operators of weight zero.
   Some of these results can be generalised for non-zero weight.

     This paper is based on the talk of one of the authors (HMK)
at the conference ``Algebraic topology and abelian functions''
in honor of Professor V.M. Buchstaber in June 2013.
  Many results of the next three sections of this article
 were already discussed 
in our paper \cite{KhBiggs1}, and  here we consider
them under a  different light. Results of the last section are new.   

\medskip

 {\bf Acknowledgment} 

We are grateful to organisers of the conference
``Algebraic topology and abelian functions''
for inviting us to this conference and encouraging us
to write this text. We are grateful also to V.Ovsienko and
   T.Voronov for useful discussions.

  \section {Pencil lifting for first order operators}

  We begin by describing the geometry of first 
order operators on algebra $\F(M)$.
(See also \cite{KhVor2}, \cite{KhVor5}.)  
 
Let $\hL$ be a first order operator 
 on the algebra of densities of weight $\delta$.
Then $\hL=\hX+F$, where $F=\hL(1)$ is a density of weight $\d$ and
$\hX$ a vector field of weight $\d$ on the extended manifold $\hM$;
     $\hX(\bs_1\cdot\bs_2)=\hX(\bs_1)\cdot\bs_2+
   \bs_1\hX(\bs_2)$. In components
 $\hX=t^\d\left(X^i(x)\p_i+X^0(x)\hw\right)$. 

  Using the canonical scalar product, \eqref{canonicalscalarproduct}, 
 one can assign to the vector field 
 $\hX$, its adjoint,  
the first order operator  
          $$
\hX^*=\left[t^\d\left(X^i(x)\p_i+X^0(x)\hw\right)\right]^*=
t^\d\left(-X^i(x)\p_i-\hw X^0-\p_iX^i(x)+(1-\d)X^0(x)\right)\,.
       $$
 Thus we come to the canonical divergence of  a 
vector field $\hX$ on 
$\hM$: 
    \begin{equation*}\label{divergence}
    {\rm div\,}\hX=-\left(\hX+\hX^*\right)=
 t^\delta\left(\p_iX^i(x)+(\delta-1)X^0(x)\right)\,.
     \end{equation*} 
In the case that $\hX=\hX(x,t)=X^i(x,t)\p_i+X^0(x,t)\hw$
is an arbitrary vector field
(weight is not fixed) then this formula reads: 
    \begin{equation}\label{divergence}
    {\rm div\,}\hX=-\left(\hX+\hX^*\right)=
 \p_iX^i(x,t)+(\hw-1)X^0(x,t)\,.
     \end{equation} 
Notice that $\hX$ is divergenceless, ${\rm div\,}\hX=0$, 
iff $\hX$ is anti-self-adjoint.

Now let us 
 focus our attention on the case when $\d=0$. 
 In this case the vector field $\hX$ projects to a vector field $\X$ 
 on the 
 manifold $M$: $\X(f):=\hX(f)$ for every function $f$ on $M$.
Similarly an  arbitrary vector field $\X$ on $M$ can be uniquely 
lifted to a vector field $\hX$ on $\hM$ such that
$\hX$ projects to $\X$ and $\hX^*=-\hX^*$ 
$\Leftrightarrow {\rm div\,}\hX=0$.
We denote such a lift of $\X$ by $\hLL_\X$:
         \begin{equation*}
X^i\p_i= \X\mapsto \hLL_\X=X^i(x)\p_i+\hw\p_i X^i(x)\,,
     \quad ({\rm div\,}\hLL_\X=0)\,.
        \end{equation*}  
  $\hLL_\X$ is the Lie derivative along the vector field $\X$: 
for every density  $\bs =s(x)|Dx|^\l$,
                   \begin{equation}\label{liederivative}
 \hLL_\X(\bs )=\left(X^i(x)\p_i+\hw\p_i X^i(x)\right)s(x)|Dx|^\l
      =\left(X^i(x)\p_is(x)+\l\p_i X^i(x) s(x)\right))|Dx|^\l\,.
                    \end{equation}
Thus we see that every vector field $\hX$ on $\hM$ of weight $\d=0$
is a sum of a 
 Lie derivative $\hLL_\X$ and a vertical vector field 
         \begin{equation*}
\hX=X^i\p_i+\hw X^0=\hLL_\X+\hw (X^0(x)-\p_iX^i(x))\,. 
        \end{equation*}  

The function $X^0(x)-\p_i X^i(x)$ is a scalar function, 
it does not change
under coordinate transformations. 
Respectively 
  every operator $\hL$ of order $\leq 1$ and
of weight $\d=0$ can be decomposed as 
  \begin{equation}\label{decompositionoffirstorderoperators}
  \hL=\hLL_\X+\hw S_1+S_2\,,
   \end{equation}
where $S_1$ and $S_2$ are scalar functions;
 $\X$ is a vector field on $M$, projection of vector field
$\hX$, $\hL=\hX+S_2$, 
and $\hX=\hLL_\X+\hw S_1$ due to \eqref{liederivative}. 
  This decomposition leads to the 
construction of canonical pencil lifting
for the space of first order operators. 

  Namely pick an arbitrary $\l$ and consider the space  
$\D^{(1)}_{\l}$ of first order operators acting on 
densities of weight $\l$ (as usual we consider only 
operators of weight $\d=0$), and 
let $L$ be an arbitrary operator in this space. 
   It defines the vector field $\X$ on $M$ 
such that $\X(f)\bs= L(f\bs )-fL(\bs)$, where $\bs$ is an arbitary
density of weight $\l$. In local coordinates
  $L=X^i(x)\p_i+F(x)$  Consider 
the Lie derivative corresponding to the vector 
field $\X$,
$\LL^{\l}_\X=\hLL_\X\big\vert_{\hw=\l}$.  The difference 
$L-\LL_\X=S$ is a zeroth order operator, a scalar 
 function $S$. Hence we arrive at the decomposition of 
the linear operator  
 $L\in \D^{(1)}_{\l}(M)$ into 
 a sum of a Lie derivative and a scalar function: 
      \begin{equation}\label{decomposition2}
 L=X^i\p_i+F(x)=\underbrace
 {\left(X^i(x)\p_i+\l\p_iX^i(x)\right)}
               _
 {\hbox{Lie derivative $\LL^{\l}_\X$}}
           +
 \underbrace{\left(F(x)-\l \p_i X^i(x)\right)}
               _
        {\hbox{scalar function $S(x)$}}\,.
      \end{equation}  
(We have supposed that $L$ is an operator of weight $\d=0$.) 

Comparing this decomposition with the decomposition 
\eqref{decompositionoffirstorderoperators} we come to the following 
canonical pencil lifting map defined on the spaces $\D^{(1)}_{\l}$:
for every $L=X^i(x)\p_i+F(x)$,
        \begin{equation}\label{givenabove}
    \hPi_\l(L)=
    \hLL_\X+(F(x)-\l\p_iX^i(x))\,.     
        \end{equation}
Using this lifting we can define a family of $\diff(M)$-equivariant liftings.
For every linear operator $L\in\D^{(1)}_\l$,
$L=X^i(x)\p_i+F(x)$, (weight $\d=0$) 
we define
        \begin{equation}\label{generalformulaforfirstorderpencilliftings}
\hPi_\l^{[p:q]}(L)=
    \underbrace{\left(X^i(x)\p_i+\hw \p_iX^i(x)\right)}_
    {\hLL_\X}+{p\hw +q\over p\l+q}
        \underbrace{
       \left(F(x)-\l\p_iX^i(x)\right)}_{S(x)
              }\,,\quad     
                    \left({q\over p}\not=-\l \right)\,,
         \end{equation}
where $[p:q]$
is a point on the projective line $\RR P^1$. 

It is evident that for all pencil liftings 
$\{\hPi^{[p:q]}_\l\}$ the condition   
\eqref{pencilpassingthroughoperator} is obeyed,
and all these pencil liftings are evidently $\diff(M)$-equivariant
since they are defined via canonical constructions.     
This construction describes all regular 
$\diff(M)$-pencil of first order operators.

The pencil lifting maps $\{\hPi^{[p:q]}_\l\}$ can be considered as a 
map on the space of parameters 
          \begin{equation*}
              \Big(
     \underbrace{\RR P^1}_{[p:q]}\times 
    \underbrace{\RR}_{\hbox{weights}}
           \Big)\backslash 
                \underbrace{\RR}_{{p\over q}\not=-\l} = \RR^2\,.
          \end{equation*}

This family possesses the following distinguished liftings:

  i) the pencil lifting $\hPi^{[p:q]}$ for $[p:q]=[0:1]$
   given above \eqref{givenabove},
         \begin{equation*}
\hPi_\l^{[0:1]}(L)=\hPi(X^i(x)\p_i+F(x))=
  \hLL_\X+S(x)\,,
         \end{equation*}
 
ii) the anti-self-conjugate pencil lifting 
  $\hPi^{[p:q]}$ for $[p:q]=[2:-1]$,
  \begin{equation*}
\hPi^\l_{[p:q]}(L)=\hPi(X^i(x)\p_i+F(x))=
  \hLL_\X+{2\hw -1\over 2\l-1}S(x)\,,\quad     
                           \l\not={1\over 2},\,\quad
   \hPi^*=-\hPi\,.
         \end{equation*}

iii)  the affine line of pencil liftings 
         \begin{equation*}
\hPi^c_\l(L)=\hPi_\l^{[p:q]}(L)=\Pi(X^i(x)\p_i+F(x))=
  \hLL_\X+S(x)+c(\hw-\l)S(x)\,,
      \end{equation*}     
 where $c$ is a point on affine line:
        \begin{equation*}
   [p:q]=\left[1:{1\over c}-\l\right]\,,\hbox{in particular}\,\,
     [p:q]=[0:1]\,\,{\rm if}\,\, c=0\,.  
         \end{equation*}

\begin{remark} The $\diff(M)$-equivariant  pencil lifting
\eqref{givenabove} is  the DLO-pencil lifting \eqref{DLOpencillifting1}
for first order operators.  
 Respectively
pencil liftings \eqref{generalformulaforfirstorderpencilliftings}
are  a special case of regular $\proj$-invariant pencil
liftings. (See in detail section \ref{projective}.) 
\end{remark}

\begin {remark} One can attempt similar analysis
for operators of arbitrary weight (for detail see \cite{Biggsdiss}) 
\end{remark}

\section{Pencil liftings of second order operators}  

 We describe the geometry of second 
order operators on the algebra of densities
  $\F(M)$.

An arbitrary  differential 
operator on the algebra of densities 
is a sum  of a self-adjoint and 
anti-self-adjoint operator:
       $$
\hD=\underbrace{\hD+\hD^*\over 2}_{\hbox {self-adjoint}}
+   \underbrace{\hD-\hD^*\over 2}_{\hbox {anti-self-adjoint}}\,.
       $$
 Let $\hD$ be an arbitrary second order 
operator on the algebra of
densities, $\F(M)$, $\hD\in \D^{(2)}(M)$ and $\hD\not\in \D^{(1)}(M)$.
Then one can see that it is 
a sum of a second order self-adjoint operator and 
first order anti-self-adjoint operator. 
   First order anti-self-adjoint operators, 
 which are related
with Lie derivatives,  were analyzed 
in the previous section. We now focus on 
 second order self-adjoint operators.
Let    
  \begin{equation}\label{secondorderoperator1}
   \hD=
    \underbrace{S^{ik}(x)\p_i\p_k+
   2 B^i(x)\p_i\hw+C(x)\hw^2}_{\hbox{second order derivatives}}+
\underbrace{D^i(x)\p_i+E(x)\hw}_{\hbox{first order derivatives}}+
          F(x)\,,\qquad (\hw=t\p_t)\,.
             \end{equation}
be a second order operator on the algebra of densities. 
 As usual we shall consider the case when the operator $\hD$ has weight 
$\delta=0$.  The  
  principal symbol of this operator is 
          \begin{equation*}\label{principalsymbol}
       \widehat S(x,t)=\begin{pmatrix}
               S^{ik}(x)& B^i(x)t\cr
                 B^k(x)t& C(x)t^2  \cr
          \end{pmatrix}\,.         
          \end{equation*}
We assign to the principal symbol $\widehat S$ the following operator
on the algebra of densities:
           \begin{equation}\label{canonicaloperator}
  \hD_{\hS}\colon\quad
   \hD_\hS (\bs)={\rm div\,}\left(\hS d\bs\right)\,,
         \end{equation}
where $\bs=\bs(x,t)$ is an arbitrary density  
which is identified with a quasi-polynomial 
on the extended manifold $\hM$,
$d\bs$ is the differential of this density,
 it is a $1$-form (covector) on the bundle space $\hM$;
  $\hS d\bs$
is the vector field on $\hM$ corresponding to the covector $d\bs$
 raised by $\hS$,
 and ${\rm div\,}$ is the canonical divergence 
\eqref{divergence}. In local coordinates 
       $$
d\bs(x,t)=\p_i\bs (x,t)dx^i+\p_t \bs(x,t)dt\,,
       $$
       $$
\hS (d\bs)=\left(
  S^{ik}(x)\p_k\bs (x,t)+
   B^i(x)t\p_t\bs (x,t)
       \right)\p_i+
       \left(
   B^k(x)t\p_k\bs (x,t)+C(x)t^2\p_t\bs (x,t)
           \right)\p_t=
       $$
       $$
=\left(
  S^{ik}(x)\p_k\bs (x,t)+
   B^i(x)\hw\bs (x,t)
       \right)\p_i+
       \left(
   B^k(x)\p_k\bs (x,t)+C(x)\hw\bs (x,t)
           \right)\hw\,,
       $$
       $$
{\rm and}\qquad  \hD_\hS (\bs)={\rm div\,}\left(\hS d\bs\right)=
       $$
           \begin{equation*}
      \label{operatortoprincipalsymbolinlocalcoordinates}
 \p_i\left(S^{ik}(x)\bs (x,t)+\hw B^i(x)\bs(x,t)\right)+
  (\hw-1)\left(B^i(x)\p_i\bs(x,t)+C(x)\left(\hw \bs (x,t)\right)\right)=
         \end{equation*} 
           \begin{equation}\label{canonicaloperator1}
    S^{ik}(x)\p_i\p_k +\p_kS^{ki}(x)\p_i+
 (2\hw-1)B^i(x)\p_i+\hw \p_kB^k(x)+
      \hw(\hw-1)C(x)\,.
            \end{equation}
It is important to note that the
canonical operator $\hD_\hS$ is self-adjoint. 
This can be checked directly. It is evident that
this operator obeys the normalisation condition $\hD_\hS(1)=0$.
Suppose now that the operator $\hD$ in 
\eqref{secondorderoperator1}
is also self-adjoint, $\hD^*=\hD$. Then the operator
$\hD-\hD_\hS$ is self-adjoint.
Both these operators have the same principal symbol.
Hence the self-adjoint operator $\hD-\hD_\hS$
is of order $\leq 1$. On the other hand an operator
of order $1$ cannot be self-adjoint. This implies that 
$\hD-\hD_\hS$ is
an operator of multiplication by a scalar function.
Thus $\hD-\hD_\hS=\hD(1)=F(x)$. 
We come to the following
\begin{proposition}\label{1}
The canonical operator $\hD=\hD_\hS$ on the algebra of densities
defined by equation \eqref{canonicaloperator}
is self-adjoint and  obeys the normalisation condition 
    \begin{equation}\label{normalisationcondition1}
            \hD(1)=0\,.
     \end{equation}
An arbitrary second order self-adjoint operator $\hD$,
$\hD^*=\hD$      
which obeys normalisation
condition \eqref{normalisationcondition1} is uniquely 
defined by its principal symbol $\hS$, $\hD=\hD_\hS$.  
\end{proposition}

\begin{remark}  Canonical divergence \eqref{divergence} 
can be considered as
the divergence with respect to 
`generalised volume form' induced by the scalar
product \eqref{canonicalscalarproduct} 
(see for details \cite{KhVor4}).  One can say
that  formula \eqref{canonicaloperator}
defines the canonical self-adjoint operator as  the
Laplace-Beltrami operator
on the extended manifold $\hM$ 
corresponding to an `upper Riemannian metric' $\hS$  
and the canonical `generalised volume form'.
\end{remark}      

   We now analyze the geometrical meaning of the components of 
the principal symbol $\hS$. Under changing 
of local coordinates
$x'=x'(x), t'=t'(x,t)=
 \det\left({\p x^{i'}(x)\over \p x^i}\right)t$ (see equations
 \eqref{changingoflocalcoordinatesinextendedmanifold})  
            $$
 \p_i={\p x^{i'}(x)\over \p x^i}\p_{i'}+\p_i\log J \hw\,,\quad
      \hw=t\p_t=t'\p_{t'}\,, \quad {\rm where}\,\,
       J=\det \left({\p x^{i'}\over \p x^i}\right)\,.
            $$
Hence components of 
the principal symbol (a symmetric contravariant tensor 
 on the extended manifold) 
transform as follows:
                    \begin{equation*}
       S^{i'k'}(x')={\p x^{i'}(x)\over \p x^i}
                    {\p x^{k'}(x)\over \p x^k}
                          S^{ik}(x(x'))\,,
                     \end{equation*}
                    \begin{equation*}
       B^{i'}(x')={\p x^{i'}(x)\over \p x^i}
                   \left(B^i(x(x'))+ 
                          S^{ik}(x(x'))\p_k\log J
                           \right)\,,
                     \end{equation*}
    \begin{equation}\label{transformationlaws}
       C'(x')=C(x(x'))+
                   2B^i(x(x'))\p_i\log J
                           +
                \p_i\log J S^{ik}(x(x'))\p_k\log J\,.
                     \end{equation}
   $S^{ik}(x)$ is a contravariant symmetric tensor 
field on the base manifold 
$M$. The operator $\hD$ defines a pencil of second 
 order operators $\{\Delta_\l\}$, $\Delta_\l=S^{ik}\p_i\p_k+\dots$,
and all these operators on $M$ have the same principal symbol $S^{ik}$.
(We suppose that the tensor field $S^{ik}(x)$ does not vanish.
If $S^{ik}\equiv 0$, then 
$\hD=(2\hw-1)\hLL_B-\hw(\hw-1)(C(x)-2\p_kB^k(x))$.

   To study the geometrical meaning of the components $B^i(x)$ and $C(x)$
of the operator \eqref{canonicaloperator} it is useful to consider 
{\it a connection} on the fibre bundle $\hM\to M$. Let  $\nabla$ be an
arbitrary connection. It defines a covariant derivative of densities
with respect to vector fields. If $\bs$ is a density of weight $\l$, 
and $\X=X^i(x)\p_i$ a vector field on $M$, then 
            \begin{equation*}
     \nabla_\X(\bs)=
  \left(X^i(x)\p_i s(x)+\l\gamma_i(x)X^i(x)s(x)\right)|Dx|^\l\,.
            \end{equation*}
For an arbitrary density  $\bs=\bs (x,t)$ 
(identified with a quasipolynomial on $\hM$) 
            \begin{equation*}
  \nabla_\X\bs(x,t)=
  \left(X^i(x)\p_i\bs(x,t)+\hw\gamma_i(x)X^i(x)\bs(x,t)\right)\,.
            \end{equation*}
The components $\gamma_i(x)$ of the connection are defined 
by the equation 
$\gamma_i(x)|Dx|=\nabla_{\p_i}|Dx|$. The 
connection $\nabla$ defines a 
lifting of vector fields from $M$ to horizontal 
 vector fields on $\hM$.
The corresponding 
connection form, $\Omega$ on $\hM$, vanishes on horizontal 
 vectors, and obeys the normalisation condition $\Omega (\hw)=1$.
In local coordinates 
               \begin{equation*}
  \X=X^i(x)\p_i\mapsto \hX=X^i(x)\p_i+\gamma_i(x)X^i(x)\hw\,,
 \quad \Omega=t^{-1}dt-\gamma_i(x)dx^i\,,\quad
 \left(\Omega(\hw)=1,\Omega(\hX)=0\right)\,.
               \end{equation*} 
In particular a volume form $\rh=\rho(x)|Dx|$ 
defines a flat connection: 
       \begin{equation}\label{flatconnection}
\nabla_\rh\colon\quad
 \Omega=d(\log\rh), \gamma_i=-\p_i\log\rho(x)\,.
              \end{equation}
Now return to the self-adjoint operator $\hD_\hS$ 
(see equations \eqref{canonicaloperator} and 
\eqref{canonicaloperator1}).
 The principal symbol $\hS$ of the canonical operator $\hD_\hS$ and
the connection form $\Omega$,  which is a 
covector in the extended manifold $\hM$,
define the vector field $\hS\Omega$ on $\hM$: 
          \begin{equation*}
   \hS\Omega=\begin{pmatrix}
               S^{ik}(x)& B^i(x)t\cr
                 B^k(x)t& C(x)t^2  \cr
          \end{pmatrix}
      \begin{pmatrix}
               -\gamma_k(x)\cr
                t^{-1}  \cr
          \end{pmatrix}=
     \left(B^i(x)-S^{ik}(x)\gamma_k(x)\right)\p_i+
      \left(C(x)-B^i(x)\gamma_i(x)\right)\hw\,.
              \end{equation*}
 Consider the projection of this vector field to the
 manifold $M$:
          \begin{equation}\label{vectorfieldvanish}
\hbox{Projection of vector field $\hS\Omega$ on $M$}=
     (B^i(x)-S^{ik}(x)\gamma_k(x))\p_i\,.
         \end{equation}
Performing these Kaluza-Klein like considerations 
(see also \cite{Kh2}) we see that $B^i(x)-S^{ik}\gamma_k(x)$
is a vector field, i.e. $B^i(x)$ is an upper connection
(see \cite{KhVor2}, \cite{KhVor4} for details).

  In the case when $B^i(x)$ belongs to the image
of $S^{ik}$ 
 one can impose the condition that 
the vector field \eqref{vectorfieldvanish} vanishes.
In the case that  $S^{ik}$ is non-degenerate then
this condition uniquely defines a connection such that
 $B^i=S^{ik}(x)\gamma_k(x)$.
 This means that
the upper connection $B^i(x)=S^{ik}(x)\gamma_k(x)$ 
is generated by a canonical connection. In this case
the vector field $\hS\Omega$
is vertical, i.e. $C(x)-\gamma_i(x)S^{ik}\gamma_k(x)$ is a scalar field.

       In the general case 
we need more considered analysis.   
 Let $\nabla$ be an arbitrary connection
 on the fibre bundle $\hM\to M$.
Consider the following second order operator defined via 
the connection $\nabla$
and the contravariant tensor field $S^{ik}(x)$:
         \begin{equation}\label{canonicaloperator2}
  \hD_{_{\S,\nabla}}(\bs)={\rm div\,}_{_{\nabla}}
      \left(\S D_{_{\nabla}} \bs\right)\,.
            \end{equation}
(Compare this operator with operator \eqref{canonicaloperator}).
Here $D_{_{\nabla}} \bs$ is a density valued $1$-form,
a weighted covector
on the manifold $M$, 
$D_{_{\nabla}} \bs=\nabla_i\bs(x,t)dx^i=
\left(\p_i \bs(x,t)+\hw \G_i(x)\bs (x,t)\right)dx^i$,
where $\G_i(x)$ are components of the connection $\nabla$; 
$\S D_{_{\nabla}}\bs=S^{ik}(x,t)\nabla_k\bs(x,t)$   
is the corresponding weighted vector field,
and  ${\rm div\,}_{\nabla}$ is divergence with respect to the
connection
$\nabla$: for an arbitrary vector field-density $\X(x,t)=X^i(x,t)\p_i$
        $$
{\rm div\,}_{_{\nabla}}\X(x,t)={\rm div\,}\hX(x,t)=
  \p_i X^i(x,t)+(\hw-1)\G_i(x)X^i(x,t)\,,
        $$
where  $\hX$ is a horizontal lifting of the vector 
 field $\X$ 
  and ${\rm div\,}$  is the canonical divergence \eqref{divergence}. 
We come to the answer:
          \begin{equation*}
\hD_{_{\S,\nabla}}=S^{ik}(x)\p_i\p_k+
    \p_r S^{ri}(x)\p_i+(2\hw-1)\G^i(x)\p_i+
     \hw\p_i\G^i(x)+\hw(\hw-1)\G^i\G_i\,,\quad (\G^i=S^{ik}\G_k)\,.
          \end{equation*}           
Comparing this expression with expression  \eqref{canonicaloperator1}
for operator \eqref{canonicaloperator} we see that 
          \begin{equation*}
  \hD_{_{\S, \nabla}}=\hD_{\hS_{\nabla}},\quad{\rm where}\quad
         \widehat S_{_{\nabla}}=
      S_{_{\nabla}}(x,t)=\begin{pmatrix}
               S^{ik}(x)& \G^i(x)t\cr
                 \G^k(x)t& \G^i(x)\G_i(x)t^2  \cr
          \end{pmatrix}\,,\qquad \G^i(x)=S^{ik}(x)\G_k(x)\,.
           \end{equation*}
(Here $\hS_{\nabla}$ can be considered as a horizontal lifting of $\S$.)
One can see that the difference between the operators 
 $\hD_\hS$ and $\hD_{\hS_{\nabla}}$ is equal to
$\hD_\hS - \hD_{\hS_{\nabla}}=$
                $$
=(2\hw-1)
      \hLL_\Y+
      \hw(\hw-1)
           \left(
  C(x)-2\G_i(x)B^i(x)+\G_i(x)\G^i(x)-
    {\rm div\,}_{_{\nabla}}\Y\right)\,,
       \, Y^i=B^i-S^{ik}\Gamma_k\,. 
                $$
We see that for the operator \eqref{canonicaloperator1},
$Y^i=B^i(x)-S^{ik}\G_k(x)$ is a vector field and
$C(x)$ is an object such that 
$C(x)-2\G_i(x)B^i(x)+\G_i(x)\G^i(x)$ is a scalar function.
This is in accordance with transformation laws 
\eqref{transformationlaws}.

\begin{remark} Above we analyzed the relation between
 second order self-adjoint operators on $\hM$ 
and their principal symbols.
In the article \cite{KhVor2} self-adjoint operators
 were studied with the use
  of a  'long bracket' which can be assigned 
to this operator. The long bracket 
$\{\quad,\quad\}_\hD$ on the space of densities is equal to 
           $$
 \{f(x,t),g(x,t)\}_\hD=\hD(fg)-f\hD g-g\hD f\,,    
           $$  
(if the normalisation condition $\hD(1)=0$ holds). 
The long bracket is
equivalent to the principal symbol, and it appears naturally when considering
the relation between odd Poisson brackets and odd second order
operators on odd symplectic supermanifolds 
(see \cite{KhVor2} and also \cite{KhVor4}). 
\end{remark}

  Recall that the problem which 
we discuss in this article is
to construct a pencil lifting of operators on densities 
of a given weight.  
  In this section we consider pencil liftings of second order operators
using the constructions above. Pick an arbitrary $\l$ and consider the space
$\D^{(2)}_\l(M)$ of operators of order $\leq 2$ 
acting on densities of weight $\l$. Consider an arbitrary operator 
         \begin{equation}\label{definedbythisequation}
\Delta=A^{ik}(x)\p_i\p_k+A^i(x)\p_i+A(x)
          \end{equation}
in $\D^{(2)}_\l(M)$. (We consider as always 
 only operators of weight $\delta=0$.) Comparing this operator with the operator
\eqref{canonicaloperator} restricted to $\hw=\l$,
 and using Proposition\ref{1}
we come to  
\begin{theorem}\label{lifting}
If $\l\not=0,{1\over 2},1$ then for an arbitrary operator 
$\Delta\in\D_\l^{(2)}(M)$
there exists a unique self-adjoint operator $\hD\in D^{(2)}(\hM)$ 
such that
\begin{itemize}
\item it obeys normalisation condition \eqref{normalisationcondition1},
$$
\hD(1)=0\,.
 $$
\item it passes through the operator $\Delta$
        $$
    \hD\big\vert_{\hw=\l}=\Delta.
        $$
\end{itemize} 
 This operator is equal to the operator $\hD_\hS$ (see equations 
\eqref{canonicaloperator} and \eqref{canonicaloperator1}) where
$\hS$ is defined as follows:
for the  operator $\Delta$ defined by equation 
\eqref{definedbythisequation}, $\hS$ is given by equations
             $$
S^{ik}(x)=A^{ik}(x)\,, B^i(x)={A^i(x)-\p_kA^{ki}(x)\over 2\l-1}\,,   
             $$
             \begin{equation}\label{components}
C(x)={A(x)-\l\p_k B^k(x)\over \l(\l-1)}=             
{A(x)\over \l(\l-1)}-{\p_k A^k(x)-\p_i\p_kA^{ki}(x)\over (\l-1)(2\l-1)}\,.
            \end{equation}
\end{theorem}

This theorem, which was first formulated in \cite{KhVor2},
 allows us to construct
a regular pencil lifting map on the spaces $\D^{(2)}_\l(M)$ 
assigning to every operator
$\Delta\in \D^{(2)}_\l(M)$ a self-adjoint operator $\hD=\hPi(\Delta)$ 
obeying the normalisation condition: if 
$\Delta=A^{ik}(x)\p_i\p_k+A^i(x)\p_i+A(x)$, then
$\hPi_\l(\Delta)$ is given by equation \eqref{canonicaloperator1}
with $S^{ik}(x), B^i(x)$ and $C(x)$ defined by equations 
\eqref{components}. We come to
       $$
 \hPi_\l(A^{ik}(x)\p_i\p_k+A^i(x)\p_i+A(x))=
 A^{ij}(x)\p_i\p_j-{2(\hw-\l)\over 2\l-1}\p_kA^{ki}(x)\p_i+
  {2\hw-1\over 2\l-1}A^i\p_i-
         $$
        \begin{equation}\label{canonicalpencilliftingmap2}
  -{\hw(\hw-\l)\over (\l-1)(2\l-1)}
   \left(\p_iA^i-\p_i\p_kA^{ki}(x)\right)+
   {\hw(\hw-1)\over \l(\l-1)}
       A(x)\,,\quad \left(\l\not=0,{1\over 2},1\right)\,.
      \end{equation} 
The theorem implies that this map is a regular pencil lifting map which is
equivariant with respect to arbitrary diffeomorphisms of $M$.

  One can consider the restriction of the image of this map to
 densities
of an arbitrary weight $\mu$:
          \begin{equation*}
        \Pi^\mu_\l=\hPi_\l\big\vert_{\hw=\mu}\,.
           \end{equation*}
The maps $\Pi^\mu_\l\colon \D^{(2)}_\l(M)\to \D^{(2)}_\mu(M)$
are equivariant with respect to the group of diffeomorphisms, 
i.e. they
are maps between $\Diff(M)$-modules 
$\D^{(2)}_\l(M)$ and $\D^{(2)}_\mu(M)$.
  In the case if $\l,\mu\not=0,{1\over 2},1$
these maps establish isomorphsims between these modules. 
Explicit expressions
for these isomorphisms immediately follow 
from equation \eqref{canonicalpencilliftingmap2}:
  If the operator $\Delta=A^{ij}(x)\p_i\p_j+A^i(x)\p_i+A(x)$ then
  its image is the operator
  $\Pi^\mu_\l(\Delta)\in \D^{(2)}_\mu(M)$  given in
  the same local coordinates by the expression
  $\Delta_\mu=B^{ij}(x)\p_i\p_j+B^i(x)\p_i+B(x)$,
    where
    \begin{equation}\label{isomorphism}
    \begin{cases}
    B^{ij}&=A^{ij}\cr
    B^i   &={2\mu-1\over 2\l-1}A^i+{2(\l-\mu)\over 2\l-1}\p_j A^{ji}\cr
    B &={\mu(\mu-1)\over \l(\l-1)}A+{\mu(\l-\mu)\over (2\l-1)(\l-1)}\left(\p_jA^j-\p_i\p_j A^{ij}\right)\cr
       \end{cases}
    \end{equation}
\begin{remark}
The theorem provides geometrical background to these isomorphisms
which first appeared in the work
\cite{DuvOvs} of Duval and Ovsienko, where they considered 
the general problem about existence of $\diff(M)$-isomorphisms beetween
spaces of second order differential operators on densities.
The authors of this work described all such isomorphisms. Later it
was observed in the article \cite{KhVor4} 
that these isomorphisms can be  described by the canonical self-adjoint
pencil constructed in \cite{KhVor2}.
\end{remark}

It is illuminating to consider the following example.

Let $\X,\Y$ be two vector fields on a manifold $M$.
Consider the second order operator $\Delta=\LL_\X\LL_\Y\in \D^{(2)}_\l$
where $\LL_\X\LL_\Y$ are Lie derivatives of densities of weight $\l$. 
We wish to 
calculate the image of this operator under the
isomorphism $\Pi^\mu_\l$.
Instead of using formulae \eqref{isomorphism} directly,
we consider a self-adjoint pencil passing 
through this operator and its restriction
to densities of weight $\mu$.   The operator
$\hLL_\X\circ\hLL_\Y$ defines a pencil of operators which passes through
the operator $\Delta=\LL_\X\circ\LL_\Y$, but this operator is not
self-adjoint.  Instead consider the operator:
                \begin{equation*}
      \hD=
  {1\over 2}\left(\hLL_\X\circ\hLL_\Y+\hLL_\Y\circ\hLL_\X\right)+
 {1\over 2} 
{2\hw-1\over 2\l-1}\left(\hLL_\X\circ\hLL_\Y-
       \hLL_\Y\circ\hLL_\X\right)\,.
               \end{equation*}
This operator is evidently self-adjoint operator since
 $\hLL_\X\circ\hLL_\Y+\hLL_\Y\circ\hLL_\X$ is self-adjoint operator
and operators $\hLL_\X\circ\hLL_\Y-\hLL_\Y\circ\hLL_\X$
and $2\hw-1$ are anti-self adjoint. On the other hand  
it defines a pencil which passes
through the operator $\Delta$: $\hD\big\vert_{\hw=\l}=
\LL_\X\circ\LL_\Y$ and normalisation condition is evidently 
obeyed. 
  The theorem implies that this formula defines the image of the
 operator
$\Delta=\LL_\X\circ\LL_\Y$ under the isomorphism $\Pi^\mu_\l$:
            $$
   \Pi^\mu_\l(\Delta)=\hD\big\vert_{\hw=\mu}=
  {1\over 2}\left(\LL_\X\circ\LL_\Y+\LL_\Y\circ\LL_\X\right)+
  {1\over 2}{2\mu-1\over 2\l-1}\left(\LL_\X\circ\LL_\Y-
  \LL_\Y\circ\LL_\X\right)=
         $$
          $$
  \LL_\X\circ\LL_\Y+
  {\mu-\l\over 2\l-1}\LL_{[\X,\Y]}\,.
            $$ 
\begin{remark}
This example first appeared in  \cite{Math2} where it was calculated
using a different technique. 
\end{remark}

\section {Pencil liftings equivariant with respect to $\sdiff(M)$}

  Now consider pencil liftings for the space
of differential operators of an arbitrary order. We would like to consider
pencil liftings $\D^{(n)}_\l\ni\Delta\to \hPi (\Delta)$ for arbitrary $n$.
In general there is no $\diff(M)$-equivariant  pencil liftings (see
 \cite{LecomteMath1}, \cite{Math2} and \cite{KhBiggs1}). We consider
pencil liftings which are equivariant with respect to 
the smaller algebra, $\sdiff(M)$, of divergenceless vector
fields. 

Let $\rh=\rho(x)|Dx|$ be a volume form on $M$. It defines a group
$\SDiff_\rh(M)$ of volume preserving diffeomorphisms and 
the corresponding Lie algebra $\sdiff_\rh(M)$.

A volume form structure on the manifold $M$ allows one
to identify spaces
of densities of different weights and operators on them
 by the following canonical isomorphisms
         \begin{equation*}
 \F_\l(M)\ni\bs(x)\to \bs'(x)=\bs(x)\rh^{\mu-\l}\in \F_\mu(M)\,,
         \end{equation*}
         \begin{equation*}
P^\mu_\l\colon\quad  \D^{(n)}_\l(M)\ni\Delta\to 
 \Delta'=P^\mu_\l(\Delta)=
\rh^{\mu-\l}\circ\Delta\circ \rh^{\l-\mu}\in\D^{(n)}_\mu(M)\,.  
        \end{equation*} 
It is useful to write down explicit formulae for isomorphisms $P^\mu_\l$
in local coordinates. If  
 $\Delta=\sum_{k=0}^n S^{i_1\dots i_k}(x)\p_{i_1}\dots\p_{k}$ is an operator
of order $\leq n$ acting on densities of weight $\l$, then
         $$
P^\mu_\l\left(\Delta\right)=
\sum_{k=0}^n \rho^{\mu-\l}(x)\circ
    S^{i_1\dots i_k}(x)\p_{i_1}\dots\p_{i_k}\circ
      \rho^{\l-\mu}(x)=
            $$
           $$
S^{i_1\dots i_k}(x)
        \left(\p_{i_1}+(\mu-\l)\G_{i_1}\right)\dots
        \left(\p_{i_k}+(\mu-\l)\G_{i_k}\right)\,,
         $$
where $\G_i=-\p_i\log\rho_i$ are components of the flat connection 
induced by the volume form $\rh=\rho(x)|Dx|$ 
(see equation \eqref{flatconnection}).
This formula follows from the observation that
        $[\p_i,\rho^\l]=\p_i\circ\rho^\l-\rho^\l\circ\p_i=-\l\G_i(x)$. 

  On a manifold with a volume form structure 
the isomorphisms $\{P^\mu_\l\}$ define a canonical pencil lifting $\hP_\l$
  on operators acting on densities of weight $\l$:
    \begin{equation}\label{canonicalpencilliftingvolume}
       \hP_\l\colon\quad \hP_\l\big\vert_{\hw=\mu}=P^\mu_\l\,,\quad
     ( \hP_\l(\Delta)=
  \rh^{\hw-\l}\circ \Delta\circ\rh^{\l-\hw}) \,.
     \end{equation}
If  
 $\Delta=\sum S^{i_1\dots i_k}(x)\p_{i_1}\dots\p_{k}$, then
         $$
\hP_\l\left(\Delta\right)=
\sum S^{i_1\dots i_k}(x)
        \left(\p_{i_1}+(\hw-\l)\G_{i_1}\right)\dots
        \left(\p_{i_k}+(\hw-\l)\G_{i_k}\right)\,.
         $$
It is evident that this canonical pencil lifting defined on operators
of an arbitrary order is equivariant with respect
to the Lie algebra $\sdiff_\rh(M)$ and the  group $\SDiff_\rh(M)$.
We now search for other pencil liftings 
which are equivariant with respect to
this algebra and to this group. 
Using these isomorphisms we may reduce 
the problem of constructing of
$\sdiff(M)$-equivariant pencil liftings to finding certain 
maps between operators on functions (densities of weight $\l=0$).

  Namely let $F=F(\Delta)$ be a linear map on the space of operators 
on functions $F\colon \D^{(n)}_0(M)\to\D^{(n)}_0(M)$. 
Then using the canonical isomorphisms  one can assign to this
map the map $\hPi^{(F)}$, defined on operators on $\F_\l(M)$
with values in operators on the algebra of densities: if
$\Delta$ is an operator on $\F_\l(M)$ then 
           \begin{equation}\label{pencilmapviamapF}
  \hPi^{(F)}_\l(\Delta)=
   \rh^\hw \circ F\Big(
\underbrace{\rh^{-\l}\circ\Delta\circ\rh^\l}_
 {\hbox{\footnotesize operator on functions}}
\Big)\circ\rh^{-\hw}
             \end{equation}
(We denote by 
$\rh^\hw\circ \Delta\circ \rh^{-\hw}$
the operator $\hD$ such that for an arbitrary weight $\mu$,
$\hD\big\vert_{\hw=\mu}=\rh^\mu\circ \Delta\circ \rh^{-\mu}$.)

The condition that the pencil $\hPi^{(F)}_\l(\Delta)$ passes through 
the operator $\Delta$, $\hPi(\Delta)\big\vert_{\hw=\l}=\Delta$  
reads that 
             \begin{equation}\label{bearinginmindcondition}
   \rh^{\l}\circ F\Big(\rh^{-\l}\circ\Delta\circ\rh^\l\Big)
               \circ\rh^{-\l}=\Delta\,.
             \end{equation}
The map $\hPi^{(F)}$ is $\sdiff(M)$-equivariant iff the map $F$
is $\sdiff(M)$-equivariant and this condition is obeyed.

  Consider the following map $F$ on operators on functions
      \begin{equation}\label{lemma1}
   F(\Delta)=\a\Delta+\beta\Delta^{*_\rh}+\gamma\Delta(1)+\delta\Delta^{*_\rh}(1)\,.
        \end{equation}
Here $\alpha,\beta,\gamma,$ and $\delta$ are constant coefficients
and $\Delta^{*_\rh}$ is the operator on functions which is adjoint
to the operator $\Delta$ with respect to the volume form $\rh$:
 $\int_M \Delta f\cdot g\rh=\int_M f\cdot (\Delta^{*_\rh}g)\rh$ 
for two arbitrary functions $f,g$ with compact support.
(Compare with the definition of the canonical adjoint (see after 
equation \eqref{canonicalscalarproduct}).)
                It is evident that this map is equivariant with respect to
all diffeomorphisms preserving the volume form and with respect to all 
divergenceless vector fields. Simple calculations show that
the corresponding map $\hPi^{(F)}(\Delta)$ 
(see equation \eqref{pencilmapviamapF})
on operators on $\F_\l(M)$  is equal to: 
          \begin{equation*}
  \hPi_\l^{(F)}(\Delta)=
a\hP_\l(\Delta)+b\left(\hP_\l(\Delta)\right)^*+
c\left(\hP_\l(\Delta)\right)(1)+d\left(\hP_\l(\Delta)\right)^*(1)\,,
       \end{equation*}               
where $\hP_\l$ is the pencil lifting 
\eqref{canonicalpencilliftingvolume} and `$^*$' denotes the canonical
adjoint.   

 Now choose an arbitrary $n$ and an arbitrary weight $\l$,
then using the last equation, we will write down a family of regular
$\sdiff_\rh(M)$-pencil liftings on the space
of operators of order $\leq n$ on densities of weight $\l$. 
Bearing in mind condition 
\eqref{bearinginmindcondition} and the fact that the 
pencil lifting is regular,
(i.e. the order of the operator $\hPi(\Delta)$ cannot be greater than $n$
if the order of the operator $\Delta$ is not greater than $n$), we come
to the fact that the coefficients $\a,\beta,\gamma,\delta$ 
are polynomials
in the weight operator $\hw$, i.e. vertical operators such that
     \begin{equation*}
      \begin{cases}
\a=A(\hw), \beta=B(\hw)\,\,\hbox{ have order $\leq 1$  in $\hw$}\,\cr
\gamma=C(\hw), \delta=D(\hw) \,\, \hbox{have order $\leq n$ in $\hw$}\cr
    A(\hw)+(-1)^n B(\hw)=1,\cr
  A(\hw)\big\vert_{\hw=\l}=1, B(\hw)\big\vert_{\hw=\l}=0\,    \cr
  C(\hw)\big\vert_{\hw=\l}=D(\hw)\big\vert_{\hw=\l}=0 \cr
 \end{cases}\qquad          
    \begin{matrix}
\hbox {the regularity condition of the pencil}\cr
                         \hbox { and the condition that}\cr
\hbox {the pencil passes through the operator}\cr
    \end{matrix}
     \end{equation*}
Thus we come to the following family of 
regular $\sdiff_\rh$-equivariant
(and $\SDiff_\rh(M)$-equivariant) pencil liftings:
$\hPi_\l$, such that for 
an arbitrary $\Delta\in \D^{(n)}_\l(M)$,    
    \begin{equation}\label{planeofliftings}
 \hPi_\l(\Delta)=
A(\hw)\hP_\l(\Delta)+B(\hw)\left(\hP_\l(\Delta)\right)^*+
C(\hw)\left(\hP_\l(\Delta)\right)(1)+
D(\hw)\left(\hP_\l(\Delta)\right)^*(1)\,,
          \end{equation}
where 
      \begin{equation*}
A(\hw)=1-b(\hw-\l)\,, B(\hw)=(-1)^nb(\hw-\l)\,,
  C(\hw)=\sum_{k=1}^n c_k(\hw-\l)^k\,,      
  D(\hw)=\sum_{k=1}^n d_k(\hw-\l)^k\,,      
\end{equation*}
$b,c_k,d_k$ are arbitrary coefficients (see also \cite{KhBiggs1}).
We come to a $2n+1$-dimensional plane of regular pencil liftings.
  The following lemma (see its proof in appendix to the article
        \cite{KhBiggs1}) 
shows that this family of pencil liftings
exhausts all regular $\sdiff(M)$-equivariant liftings (if
the manifold is of dimension greater than 2).    
\begin{lemma}
Let $F$ be a local linear map on the space $\D^{(n)}_0(M)$ of differential
operators of order $\leq n$ on functions which is equivariant 
with respect to the Lie algebra $\sdiff_\rh(M)$ of divergenceless 
vector fields. If  the manifold $M$ is connected and of
dimension greater than $2$, then $F$ is expressed by equation 
\eqref{lemma1}. 
\end{lemma}
(One should note that 
local linear maps are differential polynomials
on the coefficients of operators
(see Peetre's theorem \cite{Peetre}).)


\m

\begin{example}  In this example we consider $\sdiff(M)$-equivariant
pencil liftings of $\D^{(1)}_\l(\RR^n)$.
Let $\rh$ be a volume form in Cartesian coordinates,
$\rh=\rho(x)|Dx|$. Denote by 
$\G_i(x)=-\p_i\log\rho(x)$ the components of the flat
connection induced by the volume form (see \eqref{flatconnection}).
 If $\Delta=X^i(x)\p_i+F(x)$ is an arbitrary first order operator
acting on densities of weight $\l$, then we see that
according to equations 
 \eqref{canonicalpencilliftingvolume} 
             $$
\hP_\l(\Delta)=X^i(x)\left(\p_i+(\hw-\l)\G_i(x)\right)+F(x)\,,
             $$   
             $$
\hP^*_\l(\Delta)=-X^i(x)\p_i-\p_iX^i(x)+
    (1-\hw-\l)X^i(x)\G_i(x)+F(x)\,,
             $$   
             $$
   \hP_\l(\Delta)(1)=
\left(X^i(x)\left(\p_i+(\hw-\l)\G_i(x)\right)+F(x)\right)(1)=
    -\l X^i(x)\G_i(x)+F(x)\,,
        $$ 
             $$
   \hP^*_\l(\Delta)(1)=
    -\p_iX^i(x)+(1-\l)X^i(x)\G_i(x)+F(x)\,.
        $$ 
Thus according to \eqref{planeofliftings} we come to  
the following family of regular $\sdiff(M)$-pencil liftings: 
            $$
\hPi_\l(\Delta)=
 X^i(x)\p_i+b(\hw-\l)\p_iX^i(x)+
        (\hw-\l)(1-b(2\l-1))X^i(x)\G_i(x)+(1-2b(\hw-\l))F(x)+
         $$
         $$
       +c(\hw-\l)(F(x)-X^i(x)\G_i(x))+
       d(\hw-\l)(F(x)-\p_iX^i(x)+(1-\l)X^i(x)\G_i(x))\,,
            $$
where $b,c,d$ are constants. 
\end{example}

Now we analyse the geometrical meaning of the family 
\eqref{planeofliftings} of pencil liftings.
 One can choose a line of liftings
      \begin{equation}\label{lineofliftings}
  \{\hPi_\l^{b}\}\colon \quad
  \hPi^b_\l(\Delta)=
(1-b(\hw-\l))\hP_\l(\Delta)+(-1)^nb(\hw-\l)\left(\hP_\l(\Delta)\right)^*\,,
   b\in\RR\,.
       \end{equation}
An arbitrary lifting as in \eqref{planeofliftings}, $\hPi_\l$,
is the sum of a lifting $\hPi^{b}$ for some $b$, 
 and a map which takes values in vertical operators.
Notice that in the case if $\l\not={1\over 2}$ then the condition
      \begin{equation}\label{condition}
   b={1\over 1-2\l}
       \end{equation}
 makes the pencil lifting \eqref{lineofliftings}  
self-adjoint if $n$ is even and
anti-self-adjoint if $n$ is odd. This follows from the fact
that the vertical operators  $1-b(\hw-\l)$ and $b(\hw-\l)$ 
are the adjoints of each other 
iff condition \eqref{condition} is obeyed:
       $$
  1-b(\hw-\l)=\left(b(\hw-\l)\right)^*=b(1-\hw)-b\l\Leftrightarrow 
   b={1\over 1-2\l}\,.
       $$
Consider the pencil lifting \eqref{lineofliftings} for $b={1\over 1-2\l}$:
We will call this pencil lifting {\rm distinguished}
 and denote it $\hPi^{\rm disting}_\l$, 
          \begin{equation}\label{disting2}
         \hPi^{\rm disting}_\l(\Delta)=
   \hPi^b_\l(\Delta)\big\vert_{b={1\over 1-2\l}}=
  {\hw+\l-1\over 2\l-1}\hP_\l(\Delta)+
     (-1)^n{\l-\hw\over 2\l-1}
   \left(\hP_\l(\Delta)\right)^*\,,
          \end{equation}
               $$
   \left(\hPi^{\rm disting}_\l(\Delta)\right)^*=
    (-1)^n\hPi^{\rm disting}_\l(\Delta)\,.
           $$
Comparing the family of 
pencil liftings \eqref{planeofliftings} with the line of liftings 
\eqref{lineofliftings} and with the distinguished 
pencil lifting we see that 
that for every $\sdiff_\rh(M)$-pencil lifting 
$\hPi_\l$ the following decomposition holds:
            $$
\hbox {for all $\Delta\in \D^{(n)}_\l(M)$},
 \hPi_\l(\Delta)=(1-b(\hw-\l))\hP_\l(\Delta)+
  (-1)^nb(\hw-\l)\hP^*_\l(\Delta)+\dots=
            $$
             $$
\hPi^{b}_\l(\Delta)+\dots=
             $$
          \begin{equation}\label{decomposition}
 =\hPi_\l^{\rm disting}(\Delta)-
 k(\hw-\l)\left(
     \hP_\l(\Delta)-
 (-1)^n\hP^*_\l(\Delta)
     \right)+\dots
 \,, {\rm where}\,\, k=b-{1\over 1-2\l}\,,
         \end{equation}
  and the ellipses mean vertical operators.
It is useful at this step to introduce the following notation:
Recall that operator an $L\in \D(M)$ is vertical 
if it commutes with functions
on $M$, $L=L(t,x,\hw)$. 
We denote by $\V^{(0)}(M)$ the space of vertical 
operators on the algebra of densities.
   Inductively we define $\V^{(n)}(M)$ 
as follows: an operator $\Delta$ on $\F(M)$ 
is in $\V^{(n)}(M)$ if for an arbitrary function
$f$ on $M$, $[\Delta,f]$ is in $\V^{(n-1)}(M)$.
Roughly speaking operators
in $\V^{(n)}(M)$ are differential operators of order
$\leq n$ with coefficients which are vertical operators.
 For example consider 
the second term in
the equation \eqref{decomposition}, the operator 
$\hD=k(\hw-\l)\left(\hP_\l(\Delta)-(-1)^n\hP^*_\l(\Delta)\right)$.
 It belongs to the space $\D^n(\hM)$ of operators of order 
$\leq n$ on the algebra of densities $\F(M)$, 
but as this operator has spatial derivatives of order $\leq n-1$,
it  belongs to the space $\V^{(n-1)}(M)$.

\begin{proposition}\label{criterion}
In the case if $\l\not={1\over 2}$, then
the family \eqref{planeofliftings} of regular $\sdiff(M)$-equivariant pencil 
liftings of the space $\D^{(n)}_\l(M)$,
 contains distinguished pencil lifting \eqref{disting2}.
This lifting is self-adjoint if $n$ is even and it is anti-self-adjoint
if $n$ is odd.
The distinguished pencil lifting, which is defined on the space
$\D^{(n)}_\l(M)$
factored  by the space $\V^{(n-2)}(M)$ 
{\it does not} depend on the volume form.
In other words the terms of operator $\hPi^{\rm disting}(\Delta)$
which possess spatial derivatives of order $n$ and $n-1$
do not depend on the choice of the volume form.

   The map to the space $\D^{(n)}(\hM)$ factored by $\V^{(n-2)}(M)$
induced by the distinguished pencil lifting is
   $\diff(M)$-equivariant. 
  
  Any map to the factor space that is $\diff(M)$-equivariant and induced
   by a pencil lifting in \eqref{planeofliftings} must be induced by
the distinguished pencil lifting.

 \end{proposition}
This Proposition can be proved directly. 
Consider an example. 

\begin{example} Consider regular $\sdiff(M)$-equivariant
pencil liftings of fourth order operators.
We will perform all calculations up to the terms 
containing spatial derivatives
of order $\geq 3$, i.e. we factor out the space  $\V^{(2)}(M)$.

Let $\Delta=S^{ikmn}(x)\p_i\p_k\p_m\p_n+L^{ikm}\p_i\p_k\p_m+\dots$,
 $\Delta\in \D^{(4)}_\l(M)$, ($\l\not={1\over 2}$),
be an operator of order $\leq 4$ acting on densities
of weight $\l$ on the manifold $M$ provided with a
 volume form
$\rh=\rho(x)|Dx|$. As usual denote by $\G_i(x)=-\p_i\log\rho (x)$
the components of the flat connection corresponding to this volume form.
Then  it follows from \eqref{canonicalpencilliftingvolume} that
             $$
\hP_\l(\Delta)=S^{ikmn}(x)\p_i\p_k\p_n\p_m+
  \left (4(\hw-\l)\G_rS^{rikm}+L^{ikm}(x)\right)\p_i\p_k\p_m+\dots\,,
             $$
             $$
\hP^*_\l(\Delta)=S^{ikmn}(x)\p_i\p_k\p_n\p_m+
  \left (4(\hw+\l-1)\G_rS^{rikm}(x)+4\p_rS^{rikm}(x)-
    L^{ikm}(x)\right)\p_i\p_k\p_m+\dots\,,
             $$
and using the decomposition formulae \eqref{decomposition} we come to
           $$
\hPi_\l(\Delta)=
\underbrace{{\hw+\l-1\over 2\l-1}\hP_\l(\Delta)
                          +
            {\l-\hw\over 2\l-1}\hP^*_\l(\Delta)}
                       _
          {\hbox{operator $\hPi^{\rm disting}_\l(\Delta)$}}+
              $$
              $$
     k(\hw-\l)\left(
       \hP^*_\l(\Delta)-\hP_\l(\Delta)
                 \right)+\dots=
           $$
                        $$
                \underbrace
                     {
       S^{ikmn}(x)\p_i\p_k\p_m\p_n+
              \Big(
{\l-\hw\over 2\l-1}4\p_rS^{rikm}(x)+{2\hw-1\over 2\hl-1}L^{ikm}(x)
                \Big)\p_i\p_k\p_m
                    }
                    _
          {\hbox{operator $\hPi^{\rm disting}_\l(\Delta)$}}+
                  $$
\medskip
                  $$
     k(\hw-\l)\left(
     4\p_rS^{rikm}(x)+4(2\l-1)\G_r(x) S^{rikm}(x)-2L^{ikm}(x)
            \right)\p_i\p_k\p_m+\dots
                                 $$
We see that if $k=0$, i.e. $\hPi_\l=\hPi^{\rm disting}_\l$, then
the lifting does not depend on the volume form up to $\V^{(2)}(M)$.

\end{example}

\begin{remark}
Another proof of this proposition is  to calculate the variation
of liftings \eqref{decomposition} 
with respect to variation of the volume form,
using the formula
          $$
  \delta P_\l=(\hw-\l)[\rh^{-1}\delta\rh, P_\l]\,.
         $$
Then one can see that for an arbitrary operator $\Delta$ the variation
of the operator $\hPi^{\rm disting}(\Delta)$ lies in $\V^{(n-2)}$,
it vanishes in the factorspace (see for details \cite{KhBiggs1}).
\end{remark}

 Regular pencil liftings  equivariant with respect to   
the Lie algebra $\sdiff_\rh(M)$ can be defined on operators of
arbitrary order. The Lie algebra $\sdiff_\rh(M)$ 
is in some sense the 'closest' to the Lie algebra 
$\diff(M)$ of all vector fields.
It follows from proposition \ref{criterion} that 
any $\diff(M)$-equivariant pencil lifting must be (anti)self-adjoint
up to a simple factor. More precisely if a  $\sdiff(M)$-pencil 
lifting $\hPi$ is additionally equivariant with respect to
the Lie algebra $\diff(M)$ then the pencil lifting $\hPi$
is equal to the distinguished pencil lifting 
$\hPi^{disting}$ up to a map with values in vertical operators.
    This simple remark
becomes useful if one studies $\diff(M)$-equivariant pencil liftings
for manifolds of dimension greater than $2$, since in this case equation
\eqref{planeofliftings} describes all regular 
$\sdiff(M)$ liftings.
For example in previous sections we considered canonical pencil
liftings for first order operators and
the self-adjoint canonical lifting for second order operators.
These liftings are evidently $\diff(M)$-equivariant. 
 Considering  $\sdiff(M)$-equivariant
 distinguished pencil liftings on operators of first and second
order  and, if necessary
correcting them by vertical maps,  
one can prove that for manifolds of dimension $\geq 3$, there
are no
other $\diff(M)$-equivariant liftings except those 
considered in the previous sections
 (for details see
\cite{KhBiggs1}).

\section { The DLO-pencil lifting and 
$\proj(M)$-equivariant liftings}\label{projective}

 In this section we consider regular pencil liftings
equivariant with respect to the Lie algebra $\proj(M)$
of projective vector fields. 
 Without loss of generality we consider 
$M=\RR^d$ with Cartesian coordinates
$\{x^i\}=\{x^1,\dots,x^d\}$, the Lie algebra $\proj(\RR^d)$
is generated by vector fields
     \begin{equation*}
         \begin{cases}
         \hbox {translations\,\, $\p_i$}\cr
         \hbox {linear transformations\,\, $x^k\p_i$}\cr
         \hbox {special projective transformations\,\, 
          $x^kx^i\p_i$}\cr
         \end{cases}\quad
        i,k=1,\dots,d\,.
        \end{equation*} 
    
We will first construct strictly regular DLO-pencil lifting
(see equations \eqref{DLOpencillifting1} and
\eqref{strictlyregular} in Introduction)  
on the spaces of operators $\D^{(n)}_\l(\RR^d)$ (for arbitrary $\l$ and $n$),
then using this pencil lifting we will 
construct all regular $\proj(\RR^d)$-pencil 
liftings on $\D^{(n)}_\l(\RR^d)$.

  In the classical articles \cite{LecomteOvs1} and
\cite{DuvLecomteOvs} of Duval, Lecomte and Ovsienko  
projective and conformal equivariant quantisation
 was considered.  They constructed
full symbol maps equivariant with respect to
 the Lie algebra of projective transformations \cite{LecomteOvs1}
and full symbol maps
equivariant with respect to the Lie algebra of conformal transformations
\cite{DuvLecomteOvs}.
  These symbol maps were defined 
on operators acting on densities 
of arbitrary weight.  We consider here the projective case.
(The conformal case can be considered analogously.)
We consider the full symbol map $\s_\l$ constructed in \cite{LecomteOvs1}
which is 
 a linear map on operators
with values in contravariant symmetric tensor fields 
(polynomials in the fibre variables of $T^*M$),
 obeying the natural 
normalisation conditions (see \eqref{fullsymbolmap} in Introduction.)
For an arbitrary operator, $\Delta$,
              $$
\D^{(n)}_\l(\RR^d)\ni\Delta=\sum_{k=0}^n
       S^{i_1\dots i_k}(x)\p_{i_1}\dots\p_{i_k}\,,
              $$
              $$
 \s(\Delta)= \sum_{k=0}^n\sum_{r=0}^k c^{(k)}_r(\l)
      \p_{j_1}\dots\p_{j_r}S^{j_1\dots j_r i_1\dots i_{k-r}}
           \xi_{i_1}\dots\xi_{i_{k-r}}\,\,,
            $$
   where  the $c^{(k)}_r(\l)$ ($r=0,1,\dots,k$) are polynomials in $\l$
of order $r$ obeying the normalisation condition
                 $$
              c_0^{(k)}=1\,.
                 $$
The form of the full symbol map is dictated by the
equivariance with respect to affine transformations.
  Equivariance with respect to the special projective transformations
defines a recurrent relation amongst 
the polynomials $c_r^{(k)}(\l)$
(see for detail \cite{LecomteOvs1} and the book \cite{OvsTab}). E.g.
              $$
       \sigma(S^{ik}(x)\p_i\p_k+A^i(x)\p_i+F(x))=
   S^{ik}(x)\xi_i\xi_k-{2(\l(d+1)+1)\over d+3}\p_kS^{ki}(x)\xi_i+
            $$
          \begin{equation}\label{coefficientsforprojectivemap1}
          {\l(\l(d+1)+1)\over d+2}\p_k\p_iS^{ki}(x)+
      A^i(x)\xi_i-\l \p_iA^i(x)+F\,.
        \end{equation}
One can consider a quantisation map $Q_\l$ which is inverse 
to the full symbol map, $Q_\l=\s^{-1}_\l$,
            \begin{equation}\label{quantisationmap}
   Q_\l\left(S^{i_1\dots i_k}\xi_1\dots \xi_k\right)=
   \sum_{r=0}^k \t c^{(k)}_r(\l)\p_{j_1}\dots\p_{j_r}
     S^{j_1\dots j_r i_1\dots i_{k-r}}\p_{i_1}\p_{i_2}\dots \p_{i_{k-r}}\,,
            \end{equation}
where $\t c^{(k)}_r$ are polynomials in $\l$ which can be expressed
recurrently through the polynomials $c_k(\l)$ by the formulae 
       $\t c^{(k)}_0 c^{(k)}_0=1$ and
$\sum_{i=0}^p \t c^{(k)}_{p-i} c^{(k-p+i)}_i=0$ if
   $p>0$. We come to $\t c^{(k)}_0=1$, $\t c^{k}_1(\l)=-c^{(k)}_1(\l)$,
$\t c^{(k)}_2(\l)=-c^{(k)}_2(\l)-\t c_1^{(k)}(\l)c_1^{(k-1)}(\l)=
  -c^{(k)}_2(\l)+c_1^{(k)}(\l)c^{(k-1)}_1(\l)$, $\dots$.
In particular it follows from \eqref{coefficientsforprojectivemap1} that
               $$
       Q_\l(S^{ik}(x)\xi_i\xi_k+A^i(x)\xi_i+F(x))=
   S^{ik}(x)\p_i\p_k+a_d(\l)\p_kS^{ki}(x)\p_i+
            $$
          \begin{equation}\label{coefficientsforprojectivemap2}
         b_d(\l)\p_k\p_iS^{ki}(x)+
      A^i(x)\p_i+\l \p_iA^i(x)+F\,,
        \end{equation}
where
       \begin{equation}\label{polynomials1}
   a_d(\l)={2(\l(d+1)+1)\over d+3}\,,\quad
      b_d(\l)={\l(d+1)(\l(d+1)+1)\over (d+2)(d+3)}\,.
         \end{equation}
As was already explained in the Introduction 
(see \eqref{DLOpencillifting1})   
using the full symbol map $\s_\l$ and the quantisation map
$Q_\l$,
one can consider  $\proj(\RR^d)$-pencil lifting  $\hPi^{\rm DLO}$,
DLO-pencil lifting:
                  \begin{equation*}
      \hPi^{\rm DLO}_\l\colon \hPi^{\rm DLO}=\hQ_\hw\circ\s_\l\,,
               \end{equation*}
  where $\hQ_\hw\colon\,\,\hQ\big\vert_{\hw=\mu}=Q_\mu$,
  ($\hPi^{\rm DLO}_\l\big\vert_{\hw=\mu}=Q_\mu\circ\s_\l$). 
( It is easy to see that for every operator $\Delta\in \D^{(n)}_\l$,
$\hPi^{DLO}_\l(\Delta)$ is the differential operator on the algebra of densities
since all the coefficients $c_r^{(k)}(\l)$, $\t c_r^{(k)}(\l)$ 
are polynomials in $\l$.)
 
 One can see that the DLO-pencil lifting is a $\proj$-equivariant 
strictly regular pencil lifting which can be defined
on operators of arbitrary order. 
Our construction of this 
pencil lifting is nothing but a reformulation of 
some of the 
results due to Duval, Lecomte and Ovsienko. Our next step is 
to construct a family of
 $\proj$-equivariant pencils 'dressing' the DLO-pencil lifting
with  vertical operators. For this purpose we rewrite the
DLO-pencil lifting in a different way.      

   Consider the following canonical decomposition (with respect
to a projective structure) of the DLO-pencil lifting.
  Firstly for an arbitrary $n$ introduce the $\proj$-equivariant
 linear map $\hPi_{n,\l}$ defined on the space
 $\D^{(n)}_\l(\RR^d)$ of operators of order $\leq n$
on densities of weight $\l$:
          \begin{equation*}
\hPi_{n,\l}\colon\, \quad \D^{(n)}_\l\ni\Delta\,,\,
   \hPi_{n,\l}(\Delta)=\hQ_\hw (\spr^{(n)}(\Delta)), 
          \end{equation*} 
where $\spr^{(n)}(\Delta)$ is the principal symbol of operator $\Delta$.
If $\Delta=S^{i_1\dots i_n}\p_{i_1}\dots \p_{i_n}+\dots$
then using equation \eqref{quantisationmap} we come to 
             $$
\hPi_{n,\l}(\Delta)=
 \hQ_\hw\left(S^{i_1\dots i_n}\xi_{i_1}\dots \xi_{i_n}\right)=
   \sum \t c_r^{(n)}(\hw)\p_{j_1}\dots \p_{j_r}
   S^{j_1\dots j_r i_1\dots i_{n-r}}\p_{i_1}\dots \p_{i_{n-r}}\,.
           $$
For example  
         $$
   \hPi_{0,\l}(F(x))=F(x)\,,
        $$
         $$
   \hPi_{1,\l}\left(A^i(x)\p_i+F(x)\right)=
\hQ_\hw(A^i(x)\xi_i)=A^i(x)\p_i+
      \hw \p_i A^i(x)\,,
            $$
            $$
 \hPi_{2,\l}\left(S^{ik}(x)\p_i\p_k+A^i(x)\p_i+F(x)\right)=       
       Q_\l(S^{ik}(x)\xi_i\xi_k)=
            $$
        \begin{equation}\label{dlomap1}
   =S^{ik}(x)\p_i\p_k+a_d(\hw)\p_kS^{ki}(x)\p_i+
     b_d(\hw)\p_i\p_kS^{ik}(x)\,,
         \end{equation}
where the polynomials $a_d, b_d$ are defined by equations 
\eqref{polynomials1}.
 
  One can see that the operator 
            $
   \Delta-\hPi_{n,\l}(\Delta)\big\vert_{\hw=\l}=
       \Delta-Q_\l(\spr^{(n)}(\Delta))
           $
is an operator of order $\leq n-1$ since the operator 
$\Delta\in \Delta\in\D^{(n)}_\l$ has the order $\leq n$. 
Thus we see that  
the DLO-pencil lifting defined on operators of order $\leq n$ 
can be decomposed into a sum of $n+1oj$-equivariant maps: 
       $$
 \hPi^{DLO}_\l=\hPi_{n,\l}+\hPi_{n-1,\l}+\dots+\hPi_{0,\l}
       $$ 
and respectively an arbitrary operator $\Delta\in\D^{(n)}_\l(\RR^d)$
can be decomposed in the sum of the $n+1$ operators 
        \begin{equation}\label{decomposition1}
\Delta=\Delta_{n}+\Delta_{n-1}+\dots+\Delta_{0}\,,
       \end{equation}
such that
         $$
\hPi_\l^{\rm DLO}(\Delta)=\hPi_{n,\l}(\Delta_{n})+
   \hPi_{n-1,\l}(\Delta_{n-1})+\dots+\hPi_{0,\l}(\Delta_{0})\,,
         $$
where $\Delta_n=\hPi_{n,\l}(\Delta)\big\vert_{\hw=\l}$ and inductively
                       $
\Delta_{n-k}=\hPi_{n-k,\l}\left(
\Delta-\Delta_{n}-\dots-\Delta_{n-k+1}
                  \right)\big\vert_{\hw=\l}
                        $, and
                 $$
\hPi^{\rm DLO}_\l(\Delta)=\hPi_{n,\l}(\Delta_n)+
 \hPi_{n-1,\l}(\Delta_{n-1})+\dots+
  \hPi_{1,\l}(\Delta_1)+\hPi_{0,\l}(\Delta_0)\,.
                 $$
(In fact $\hPi^{\rm DLO}_\l(\Delta_k)=\hPi_{k,\l}(\Delta_k)$).

All the operators of the pencil $\hPi^{\rm DLO}_\l(\Delta_n)$ 
including
the operator $\Delta_n$ have the same
full symbol  which is equal to their principal symbol:
If in decomposition \eqref{decomposition1},
 principal symbol of $\Delta_k$ is equal to 
$S^{i_1\dots i_k}\xi_{i_1}\dots \xi_{i_k}$ then 
for an arbitrary weight $\mu$
              $$
   \sigma_\mu\left(\hPi_{n,\l}(\Delta_n)\big\vert_{\hw=\mu}\right)
     =S^{i_1\dots i_k}\xi_{i_1}\dots \xi_{i_k}\,.
             $$
 One can see that all operators of this pencil and their
 symbol are eigenvectors of the Casimir operator of  the Lie algebra 
$\proj(\RR^d)$ with the same eigenvalue.

\begin{example}\label{secondorderoperatordecomposition} 
Consider second order operator 
$\Delta=S^{ik}(x)\p_i\p_k+A^i(x)\p_i+F(x)\in \D^{(2)}_\l(\RR^d)$
acting on densities of weight $\l$.

The operators $\hPi_{2,\l},\hPi_{1,\l}$ and $\hPi_{0,\l}$ on the algebra
of densities (operator pencils) are defined
by equation \eqref{dlomap1}.
 
 Consider the decomposition $\Delta=\Delta_2+\Delta_1+\Delta_0$.
Due to the formulae above and equation 
\eqref{coefficientsforprojectivemap2} we have 
       $$
\Delta_2= 
 S^{ik}(x)\p_i\p_k+a_d(\l)\p_kS^{ki}(x)\p_i+b_d(\l)\p_k\p_iS^{ki}(x)\,,
        $$
      $$
\Delta_1=\left(A^i(x)-a_d(\l)\p_kS^{ki}(x)\right)\p_i+
   \l\left(\p_iA^i(x)-a_d(\l)\p_i\p_kS^{ki}(x)\right)\,,
      $$
         $$
\Delta_0=F(x)-\l\left(\p_iA^i(x)-a_d(\l)\p_i\p_kS^{ki}(x)\right)
    -b_d(\l)\p_k\p_iS^{ki}(x)
        $$
and respectively
       $$
\hPi^{\rm DLO}_\l(\Delta)=\hD_2+\hD_1+\hD_0\,,
       $$
where
        $$
\hD_2=\hPi^{\rm DLO}_\l(\Delta_2)=
 S^{ik}(x)\p_i\p_k+a_d(\hw)\p_kS^{ki}\p_i+b_d(\hw)\p_k\p_iS^{ki}(x)\,,
        $$
       $$
\hD_1=\hPi^{\rm DLO}_{\l}(\Delta_1)= 
\left(A^i(x)-a_d(\l)\p_kS^{ki}(x)\right)\p_i+
   \hw\left(\p_iA^i(x)-a_d(\l)\p_i\p_kS^{ki}(x)\right)\,,
    $$
     \begin{equation}\label{secondorderdecomposition}
\hD_0=\hPi^{\rm DLO}_{\l}(\Delta_0)=\Delta_0=
     F(x)-\l\p_i A^i(x)+(\l a_d(\l)-b_d(\l))\p_i\p_k S^{ki}(x)\,.
       \end{equation}
 (Polynomials $a_d,b_d$ are defined by equations
\eqref{polynomials1}.)
\end{example}

Now using the decomposition above
we consider the following dressing of the DLO-pencil lifting.

Pick an arbitrary $n,\l$ and take an arbitrary  lower 
 triangular $(n+1\times n+1)$ matrix of numbers
$||a_{ij}||$, where $i,j=0,1,2,3,\dots$ such that in every row
at least one of the elements is not equal to zero
         \begin{equation}
||a_{ij}||\colon\quad
     \begin{cases}
      a_{ij}=0, \, {\rm if}\,\, j>i\,,\cr
   \hbox{in every row there exists a non-zero element}\,\cr
       \end{cases} \,. 
        \end{equation}
  Now take the DLO-pencil lifting of the space $\D^{(n)}_\l(\RR^d)$
 and assign to every matrix obeying
the conditions above the pencil:
  $\hPi^{||a_{ij}||}_\l$ such that
              $$
\hPi^{||a_{ij}||}_\l(\Delta)=
 \sum_{i=0}^n {P_i(\hw)\over P_i(\l)}
          \hPi^{\rm DLO}_\l(\Delta_{n-i})=
                    $$
             \begin{equation}\label{triangular}
    \hPi^{\rm DLO}(\Delta_n)+
{a_{10}+a_{11}\hw\over a_{10}+a_{11}\l}\hPi^{\rm DLO}_\l(\Delta_{n-1})+
 {a_{20}+a_{21}\hw+a_{22}\hw^2\over a_{20}+a_{21}\l+a_{22}\l^2}
\hPi^{\rm DLO}_\l(\Delta_{n-2})+\dots\,,
            \end{equation}
where $\Delta=\Delta_n+\Delta_{n-1}+\dots+\Delta_1+\Delta_0$
is decomposition \eqref{decomposition1} 
 and the polynomials $P_i(\l)$ are defined by the matrix $||a_{ij}||$:
          $P_i(\l)=\sum a_{ij}\l^j$.

 One can see that all these pencil liftings are regular pencil liftings.
If for two  matrices $||a_{ik}||$ and $||a'_{ik}||$ obeying conditions
\eqref{triangular} all rows are proportional then the corresponding pencils
coincide.  This formula describes the family
 of pencil liftings parameterised by 
            $$
\RR P^1\times \RR P^2\times\dots \RR P^n\,.
            $$ 
On the other hand let $\hPi_\l$ be an 
arbitrary regular $\proj(\RR^d)$-equivariant pencil lifting
on the space $\D^{(n)}_\l(M)$. It follows from 
decomposition \eqref{decomposition1} that for every $k$,
($k=1,\dots,n$, $\Delta=\Delta_n+\Delta_{n-1}+\dots+\Delta_0$),  
        $
\hPi_\l(\Delta_k)\big\vert_{\hw=\mu}=
F_{n-k}(\mu)\Pi_k(\Delta_k)$, where $F_{n-k}$ is a polinomial such $F_{n-k}(\l)=1$.
The conditions of regularity of the pencil dictate that 
$F_{n-k}$ are polynomials of order $\leq n-k$,
$F_{n-k}(\hw)=\sum_{r=0}^k a_{nr}\hw^r$.  Thus we come to the fact
that $\hPi_\l=\hPi^{||a_{ik}||}$. We come to the proposition

\begin{proposition}
For arbitrary $n,\l$
equation \eqref{triangular} describes the  family
  of regular $\proj$-equivariant pencil 
liftings of operators $\D^{(n)}_\l(M)$. These are all regular 
$\proj$-equivariant pencil liftings of $\D^{(n)}_\l(M)$. 

\end{proposition}

One can see from decomposition \eqref{decomposition1}
and  uniqueness arguments that 
the DLO-pencil lifting is self-adjoint on operators $\Delta_{2k}$
  and anti-self-adjoint on operators $\Delta_{2k+1}$. 
On the other hand the polynomial $P(\hw)$, a vertical operator,
is self-adjoint iff 
$P_k(\hw)$ is an even polynomial of the variable $\hw'=\hw-{1\over 2}$
and it is anti-self-adjoint if it is an odd polynomial of
$\hw'=\hw-{1\over 2}$ since $\hw'=\hw-{1\over 2}$ is 
an anti-self-adjoint operator. Thus we come to a description of
self-adjoint and anti-self-adjoint projective liftings.

\begin{proposition} Let $\hPi$ be a regular $\proj$-equivariant
lifting on the space $\D^{(n)}_\l(M)$.
Then this lifting is self adjoint, i.e. it takes values
in self-adjoint operator pencils,
 if $n$ is even and in the equation \eqref{triangular}
 for all the polynomials
$\{P_{i}(\hw)\}$ 
     \begin{equation}\label{condition28}
     P_i(\hw)^*=(-1)^iP_i(\hw)\,,  
  \end{equation}
i.e. polynomials $P_{2i}$ are self-adjoint
      \begin{equation*}\label{condition30}
     P_{2i}(\hw)=\sum c_k\left(\hw-{1\over 2}\right)^{2k}\,,
  \end{equation*}
and polynomials $P_{2i+1}$ are anti-self-adjoint
      \begin{equation}\label{condition30}
     P_{2i+1}(\hw)=\sum c_k\left(\hw-{1\over 2}\right)^{2k+1}\,.
  \end{equation}

  Respectively the regular $\proj$-equivariant
pencil lifting
is anti-self-adjoint if $n$ is odd and 
if all the conditions above for polynomials are fullfilled.

\end{proposition}

It is very illuminating to consider the following example of 
self-adjoint
regular $\proj$-liftings for second order operators. 

  An anti-self-adjoint polynomial of order $\leq 1$
 is proportional to $2\hw-1$ 
and a self-adjoint polynomial of order $\leq 2$ is proportional
to
     $$
a+b\left(\hw-{1\over 2}\right)^2=a+{b\over 4}+b\hw(\hw-1)\,.  
      $$ 
Hence an arbitrary self-adjoint regular $\proj$-equivariant pencil lifting
of second order operators on densities of weight $\l$ is of the form  
      $$
  \hPi^{\rm DLO}_{\l}(\Delta_2)+
   {2\hw-1\over 2\l-1}\hPi^{\rm DLO}_{\l}(\Delta_1)
     +{p+q\hw(\hw-1)\over p+q\l(\l-1)}
        \Delta_0\,,
       $$
where $\Delta=\Delta_2+\Delta_1+\Delta_0$ is 
decomposition \eqref{decomposition1}.
We come to a projective line of self-adjoint liftings
on densities of weight $\l$ 
($p\not=0$ or $q\not=0$).

Take an arbitrary second order operator acting 
on densities of weight $\l$
      \begin{equation*}
\Delta=S^{ik}(x)\p_i\p_k+A^i(x)\p_i+F(x)\,.
      \end{equation*}
According to equations 
\eqref{secondorderdecomposition} in  
 example \ref{secondorderoperatordecomposition}
  we have 
      $$
\hPi^{([p:q])}(\Delta)= 
       \left(
      \underbrace{
 S^{ik}(x)\p_i\p_k+a_d(\hw)\p_kS^{ki}(x)\p_i+
b_d(\hw)\p_i\p_kS^{ki}(x)}_{\hD_2=\hPi^{\rm DLO}_\l(\Delta_2)}
      \right)+
        $$
        $$
     + {2\hw-1\over 2\l-1}
         \left(
      \underbrace{
     \left(
  A^i(x)-a_d(\l)\p_kS^{ki}(x)
      \right)\p_i+
    \hw\left(
  \p_iA^i(x)-a_d(\l)\p_i\p_k S^{ki}(x)
      \right)}_{\hD_1=\hPi^{DLO}_\l(\Delta_1)}
       \right)+
         $$
           \begin{equation}\label{projectiveselfadjoint1}
     +{p+q\hw(\hw-1)\over p+q\l(\l-1)}
      \left(
  \underbrace{
 F(x)-\l \p_iA^i(x)+(\l a_d(\l)-b_d(\l))
  \p_k\p_i S^{ki}(x)}_{\Delta_0}\right)\,.
        \end{equation}

   We now compare pencils in this family 
with  the canonical self-adjoint lifting given by 
theorem \ref{1}. Recall that for the operator
$\Delta=S^{ik}(x)\p_i+A^{i}(x)\p_i+F(x)$ the self-adjoint 
canonical lifting
is equal to 
                \begin{equation}\label{canonicalselfadjoint2}
 \hD=S^{ik}(x)+\p_k S^{ik}(x)\p_i+(2\hw-1)B^i(x)\p_i+
\hw\p_iB^i(x)+\hw(\hw-1)C(x)\,,
                \end{equation}
where
         $$
B^i(x)={A^i(x)-\p_kS^{ik}(x)\over 2\l-1}\,,
     C(x)={F(x)\over \l(\l-1)}-
 {\p_iA^i(x)-\p_i\p_k S^{ki}(x)\over (\l-1)(2\l-1)}\,.
        $$
(See equations \eqref{canonicaloperator1} and \eqref{components}.
We suppose that $\l\not=0,{1\over 2},1$.)

 Note that in the case when parameter for lifting 
\eqref{projectiveselfadjoint1} is equal to $[p:q]=[0:1]$.
In this case the lifted operator
$\hD=\hPi^{([0:1])}(\Delta)$ 
obeys the normalisation condition $\hD(1)=0$, 
hence due to theorem \ref{lifting} it is equal
to self-adjoint lifting \eqref{canonicalselfadjoint2}.

 Now comparing an
arbitrary self-adjoint operator $\hPi^{([p:q])}(\Delta)$ from the family
\eqref{projectiveselfadjoint1} with  the 
 canonical self-adjoint operator  $\hD=\hPi^{([0:1])}$
we see that their difference is a vertical operator. Hence we come to
decomposition:
               $$
 \hPi^{([p:q])}(\Delta)=\hD+
      \left(
  {p+q\hw(\hw-1)\over p+q\l(\l-1)}-
  {\hw(\hw-1)\over \l(\l-1)}
           \right)
   \left(\underbrace{
 F(x)-\l \p_iA^i(x)+
 (\l a_d(\l)-b_d(\l))\p_i\p_kS^{ki}(x)}_{{\cal S}}\right)\,. 
               $$
The expression 
          $$
{\cal S}={\cal S}_\Delta=F(x)-\l \p_iA^i(x)+(\l a_d(\l)
          -b_d(\l))\p_i\p_kS^{ki}(x)\,,
          $$
due to its construction 
transforms as a scalar under arbitrary projective transformations.
From the uniqueness of the canonical self-adjoint pencil it follows that
under arbitrary diffeomorphisms $\cal S$ transforms
non-trivially. We come to a cocycle which is a multidimensional 
analogue of the usual Schwarzian (see the book \cite{OvsTab}).

\begin{remark}
It is instructive to compare calculations in this example
with those in $\cal x$ 7.1 of the book \cite{OvsTab}.
\end{remark}

\medskip

  Let us finish with a remark. If a pencil lifting factors through
 the space of symbols then it is induced by a full symbol calculus
(see the Introduction). In this article the only such pencil lifting
is the DLO-lifting. The projective quantisation map inducing this lifting
is unique. However there are many regular $\proj$-equivariant liftings
of differential operators, these are maps between 
non-commutative algebras.

\end{document}